\documentclass[reqno,11pt,a4paper]{article}
\usepackage{amsmath,amsthm,amssymb,color,graphicx,bm,mathtools}
\usepackage{algorithm}
\usepackage{algorithmic}

\newtheorem{theorem}{Theorem}[section]

\newtheorem{definition}[theorem]{Definition}

\theoremstyle{remark}
\newtheorem{remark}[theorem]{Remark}
\newtheorem{example}[theorem]{Example}

\numberwithin{equation}{section}
\numberwithin{theorem}{section}

\renewcommand{\phi}{\varphi}

\newcommand{\nada}[1]{}

\newcommand{\R}{\mathbb R}

\newcommand{\beq}{\begin{equation}}
\newcommand{\eeq}{\end{equation}}

\let\div\undefined
\DeclareMathOperator{\div}{div}
\DeclarePairedDelimiter{\norm}{\lVert}{\rVert}

\definecolor{darkgreen}{rgb}{0,0.55,0}

\newcommand{\res}{\mathop{\hbox{\vrule height 7pt width .5pt depth
               0pt\vrule height .5pt width 6pt depth 0pt}}\nolimits}
\newfont{\indic}{bbmss12}

\definecolor{light}{gray}{.97}

\title{A convex approach to the Gilbert--Steiner problem}
\author{M. Bonafini\thanks{Dipartimento di  Matematica, Universit\`a di Trento,
Italy, e-mail: mauro.bonafini@unitn.it},
\'E. Oudet\thanks{Laboratoire Jean Kuntzmann, Universit\'e de Grenoble Alpes,
 France, e-mail: edouard.oudet@imag.fr}}
\date{\today}

\begin{document}

\maketitle

\begin{abstract}

We describe a convex relaxation for the Gilbert--Steiner problem both in $\R^d$ and on manifolds, extending the framework proposed in \cite{BoOrOu}, and we discuss its sharpness by means of calibration type arguments. The minimization of the resulting problem is then tackled numerically and we present results for an extensive set of examples. In particular we are able to address the Steiner tree problem on surfaces.
\end{abstract}


\section{Introduction}



In the Steiner tree problem, at least in its classical Euclidean version, we are
given $N$ distinct points $P_1, \dots, P_N$ in $\R^d$ and we have to find the
shortest connected graph containing the points $P_i$. From an abstract point of
view this amounts to find a graph solving the variational problem
\[
\inf \{ \mathcal{H}^1(L),\; L \text{ connected},\; L \supset \{P_1, \dots, P_N\} \},\leqno{(\text{STP})}
\]
where $\mathcal H^1$ denotes the one dimensional Hausdorff measure in $\R^d$. An
optimal (not necessarily unique) graph $L$ always exists and, by minimality, $L$
is indeed a tree. Every optimal tree can be described as a union of segments
connecting the endpoints and possibly meeting at $120^\circ$ in at most $N-2$
further branch points, called Steiner points.

On the other hand, the (single sink) Gilbert--Steiner problem \cite{Gi} consists in finding a network $L$ along which to flow unit
masses located at the sources $P_1, \dots, P_{N-1}$ to the unique target point
$P_N$. Such a network $L$ can be viewed as $L = \cup_{i=1}^{N-1}\lambda_i$, with
$\lambda_i$ a path connecting $P_i$ to $P_N$, corresponding to the trajectory of
the particle located at $P_i$. To favour branching, one is led to optimize a
cost which is a sublinear (concave) function of the mass
density $\theta(x) = \sum_{i=1}^{N-1} \mathbf{1}_{\lambda_i}(x)$: i.e., for
$0\leq\alpha\le 1$, find
$$
\inf \left\{ E^\alpha(L) = \int_L |\theta(x)|^\alpha d{\mathcal H}^1(x) \right\}.
\leqno{(I_\alpha)}
$$
Problem $(I_\alpha)$ can be seen as a particular instance of an $\alpha$-irrigation problem
\cite{BeCaMo, Xia} involving the irrigation of the atomic measures $\sum_{i=1}^{N-1} \delta_{P_i}$ and $(N-1)\delta_{P_N}$, and we notice that $(I_1)$ corresponds to the Monge optimal transport problem, while
$(I_0)$ corresponds to (STP) (the energy to be optimized reduces to the length of $L$). As for (STP) a solution to $(I_\alpha)$ is known to exist and any optimal network $L$ turns out to be a tree \cite{BeCaMo}.

The Steiner tree problem is known to be computationally hard (even NP complete
in certain cases \cite{Karp}), nonetheless in $\R^2$ and $\R^3$ we have
efficient algorithms which allow us to obtain explicit solutions (see for
instance \cite{warme2001geosteiner, MINLP}), while a comprehensive survey on
PTAS algorithms for (STP) can be found in \cite{Ar, Ar1}. However, the general
applicability of these schemes restricts somehow to the Steiner tree case. For
this reason we stick here with a more abstract variational point of view, which
allows us to treat in a unified way the Steiner and Gilbert--Steiner problems.


Many different variational approximations for (STP) and/or $(I_\alpha)$ have
been proposed, starting form the simple situation where the points $P_i$ lie on
the boundary of a convex set: in this case (STP) is known to be an instance of
an optimal partition problem \cite{AmBr, AmBr2}. More recently several authors
treated these problems in the spirit of $\Gamma$-convergence using approximating
functionals modelled on Modica--Mortola or Ambrosio--Tortorelli type energies,
initially focusing mainly on the two dimensional case \cite{OuSa, BoLeSa,
	ChMeFe}, lately extending the same ideas also to higher dimensions
\cite{ChMeFe2, BoBrLe}.

Within this sole we introduce in \cite{BoOrOu} a $\Gamma$-convergence type
result in the planar case and at the same time we propose a convex framework
for the Steiner and Gilbert--Steiner problem. The approach moves from the work
of Marchese and Massaccesi \cite{MaMa} and considers ideas from
\cite{ChCrPo} in order to obtain a convex relaxation of the energy we are
dealing with. The aim of this paper is then to provide an extensive numerical
investigation of the relaxation proposed in \cite{BoOrOu}, adapting it to the
treatment of more general Gilbert--Steiner problems (with multiple sources/sinks) and addressing its validity
and applicability to problems defined on manifolds. In contrast to classical
$\Gamma$-convergence type approaches, which may numerically end up in local
minima (unless carefully taking initial guesses), this convex formulation is
able to identify (in many cases) convex combinations of optimal networks,
allowing us to have an idea of their structure. Furthermore, up to our knowledge, this is the very first formulation leading to a numerical approximation of the Steiner tree problem on manifolds.



The paper is organized as follows. In Section \ref{sec:convexrelaxation} we
review the convex framework presented in \cite{BoOrOu} for the
$\alpha$-irrigation problem $(I_\alpha)$ and extend it to the treatment of more general situations with multiple sources/sinks, both in $\R^d$ and on manifolds. In Section \ref{sec:graph} we see how the formulation simplifies for a network (STP) on graphs, with the relevant energy reducing to the norm introduced in \cite{MaMa}.
We then proceed in Section \ref{sec:algo} to describe our algorithmic scheme for the
minimization of the proposed energy functional in the Euclidean setting and we present in Section \ref{sec:results} various results for (STP) and $\alpha$-irrigation problems in two and three dimensions. In Section \ref{sec:surfaces} we eventually detail our algorithmic approach on surfaces and present some results obtained on spheres, tori and other surfaces with boundaries.

\section{Convex relaxation for irrigation type problems}\label{sec:convexrelaxation}
In this section we first review the convex framework introduced in \cite{BoOrOu}
for the $\alpha$-irrigation problem $(I_\alpha)$ and
then discuss how this same formulation can be extended to address more general Gilbert--Steiner problems with multiple sources/sinks in $\R^d$ or even on manifolds.

\subsection{The Euclidean Gilbert--Steiner problem}
Fix a set of $N$ distinct points $A=\{P_1,\dots,P_N\}\subset\R^d$, $d\ge2$.
A candidate minimizer for $(I_\alpha)$ is given as a family of simple rectifiable curves
$(\gamma_i)_{i=1}^{N-1}$, each one connecting $P_i$ to $P_N$. For optimality
reasons we can choose these curves so that the resulting network $L = \cup_i
\lambda_i$ contains no cycles (see Lemma 2.1 in \cite{MaMa}), restricting this way
the set of possible minimizers to the set of (connected) \textit{acyclic
graphs} $L$ that can be described as
\[
L=\bigcup_{i=1}^{N-1} \lambda_i, \quad\text{s.t.} \quad
\begin{aligned}
& \text{$\cdot$ }\lambda_i \text{ is a simple rectifiable curve connecting } P_i
\text{ to } P_N, \\
& \text{$\cdot$ }\text{each $\lambda_i$ is oriented by an $\mathcal{H}^1$-measurable
unit vector field $\tau_i$}, \\
& \text{$\cdot$ }\text{$\tau_i(x)=\tau_j(x)$ for $\mathcal H^1$-a.e. $x\in\,
\lambda_i\cap\lambda_j$},
\end{aligned}
\]
where the last condition requires the $N-1$ pieces composing $L$ to share the
 same orientation on intersections. Let us call $\mathcal{G}(A)$ the set of
 acyclic graphs $L$ having such a representation. Hence, we can reduce ourself
 to consider
$$
\inf \left\{ \int_L |\theta(x)|^\alpha d{\mathcal H}^1, \quad L \in \mathcal{G}(A),
 \;\;\theta(x) = \sum_{i=1}^{N-1} \mathbf{1}_{\lambda_i}(x) \right\}.
$$

To each $L\in \mathcal{G}(A)$ we now associate a measure taking values in
$\R^{d\times (N-1)}$ as follows: identify the curves $\lambda_i$ with the vector
 measures $\Lambda_i=\tau_i\cdot \mathcal{H}^1\res\lambda_i$, and consider the
 rank one tensor valued measure $\Lambda = (\Lambda_1, \dots, \Lambda_{N-1})$,
 which can be written as $\Lambda=\tau\otimes g\cdot \mathcal H^1\res L$, with
\begin{itemize}
\item $\tau \colon \R^d \to \R^d$ a unit vector field providing a global
 orientation for $L$, satisfying $\text{spt}\, \tau=L$ and $\tau=\tau_i \;
  \mathcal H^1\text{-a.e. on } \lambda_i$,
\item $g \colon \R^d \to \R^{N-1}$ a multiplicity function whose entries satisfy
$g_i\cdot\mathcal H^1\res L=\mathcal H^1\res\lambda_i$.
\end{itemize}
Observe that $g_i\in\{0,1\}$ a.e. for any $1\le i\le N-1$ (in particular
$g_i(x) = 1$ if $x \in \lambda_i$), and by construction the measures $\Lambda_i$
 verify
\begin{equation}\label{eq:solenoidal}
\div \Lambda_i=\delta_{P_i}-\delta_{P_N}.
\end{equation}

\begin{definition}\label{def:tensor} Given any graph $L \in \mathcal{G}(A)$, we
     call the above constructed measure $\Lambda=\tau\otimes g\,\cdot\,\mathcal
     H^1\res L$ the canonical (rank one) tensor valued measure representation of
     the acyclic graph $L$ and denote the set of such measures as $\mathcal{L}(A)$.
\end{definition}
Let us define on the space of matrix valued Radon measures $\mathcal{M}(\R^d; \R^{d\times(N-1)})$ the functional
\[
\mathcal F^\alpha(\Lambda) =
\left\{
\begin{aligned}
&\int_{\R^d} ||g||_{1/\alpha} \,d\mathcal H^1\res L &\quad& \text{if } \Lambda=\tau\otimes g \cdot \mathcal H^1\res L \in \mathcal{L}(A) \\
& +\infty &\quad& \text{otherwise}
\end{aligned}
\right.
\]
where we assume $1/0 = \infty$. When $\Lambda=\tau\otimes g \cdot \mathcal H^1\res L \in \mathcal{L}(A)$, since by construction $g_i \in \{0,1\}$ on $L$
and $g_i(x) = 1$ whenever $x \in \lambda_i$, one immediately gets
\[
\mathcal F^\alpha(\Lambda) = \int_{L} \left( \sum_{i=1}^{N-1} g_i(x)^{1/\alpha}
\right)^\alpha \,d\mathcal H^1 = \int_{L} \left( \sum_{i=1}^{N-1} g_i(x)
\right)^\alpha \,d\mathcal H^1 = \int_{L} \left( \sum_{i=1}^{N-1}
 \mathbf{1}_{\lambda_i}(x) \right)^\alpha \,d\mathcal H^1,
\]
which is exactly the cost $E^\alpha$ associated to $L$ in $(I_\alpha)$.
We recognize that minimizing $\mathcal F^\alpha$ among measures $\Lambda \in
\mathcal{L}(A)$ corresponds to minimize $E^\alpha$ among graphs
$L \in \mathcal{G}(A)$, and thus solves $(I_\alpha)$ in $\R^d$.

This reformulation of $(I_\alpha)$ involves the minimization of a convex energy,
namely $\mathcal F^\alpha$, but the problem is still non convex due to the non convexity of $\mathcal{L}(A)$ (the domain of definition of $\mathcal{F}^\alpha$). In view of a convex formulation the optimal choice would be to
consider the convex envelope $(\mathcal{F}^\alpha)^{**}$ of the energy, but
such an object (up to our knowledge) has no explicit
representation. Hence, following \cite{ChCrPo}, we instead look for a ``local'' convex envelope of the form
\begin{equation}\label{eq:genralRa}
\mathcal{R}^\alpha(\Lambda) = \int_{\R^d} \Psi_\alpha(\Lambda)
\end{equation}
with $\Psi_\alpha \colon \R^{d\times (N-1)} \to [0, +\infty)$ a $1$-homogeneous, convex, continuous function such that $\mathcal{R}^\alpha(\Lambda) = \mathcal{F}^\alpha(\Lambda)$ whenever $\Lambda \in \mathcal{L}(A)$. The integral in \eqref{eq:genralRa}, as outlined in \cite{BoVa}, can be defined as
\begin{equation}\label{eq:defpsilambda}
\begin{aligned}
&\int_{\R^d} \Psi_\alpha(\Lambda) = \int_{\R^d} \Psi_\alpha\left( \frac{d\Lambda_a}{d\mathcal{L}^d} \right)\,dx + \int_{\R^d} \Psi_\alpha \left( \frac{d\Lambda_s}{d|\Lambda_s|} \right)\,d|\Lambda_s| \\
&= \sup_{\phi \in C^\infty_c \left(\R^d; \R^{d\times (N-1)}\right)} \left\{ \sum_{i=1}^{N-1} \int_{\R^d} \phi_i\,d\Lambda_i - \int_{\R^d} \Psi_\alpha^*(\phi) \,dx,\quad \Psi_\alpha^*(\phi) \in L^1(\R^d) \right\}
\end{aligned}
\end{equation}
where $\Lambda = \Lambda_a + \Lambda_s$ is the Lebesgue decomposition of $\Lambda$ w.r.t. the $d$-dimensional Lebesgue measure $\mathcal{L}^d$, $|\Lambda_s|$ is the total variation of $\Lambda_s$, $\phi_i$ are the columns of the function $\phi(x) = (\phi_1(x),\dots,\phi_{N-1}(x))$ and $\Psi_\alpha^*$ is the Legendre-Fenchel conjugate of $\Psi_\alpha$ on $\R^{d\times(N-1)}$: for $p=(p_1,\dots,p_{N-1}) \in \R^{d\times(N-1)}$ and $q = (q_1,\dots,q_{N-1}) \in \R^{d\times(N-1)}$ we have
\[
\Psi_\alpha^*(q) = \sup_p \left[\left\langle q, p \right\rangle - \Psi_\alpha(p)\right] = \sup_p \left[\sum_{i=1}^{N-1} q_i \cdot p_i - \Psi_\alpha(p)\right].
\]
We immediately see that the evaluation of $\mathcal{R}^\alpha$ on any $\Lambda \in \mathcal{L}(A)$, i.e. $\Lambda = \tau \otimes g \cdot \mathcal{H}^1\res L$ with $||\tau||_2=1$ and $g_i \in \{0,1\}$, only involves the singular part of the decomposition, so that
\[
\mathcal{R}^\alpha(\Lambda) = \int_{\R^d} \Psi_\alpha(\tau \otimes g) d\mathcal{H}^1\res L.
\]
Since we require $\mathcal{R}^\alpha(\Lambda) = \mathcal{F}^\alpha(\Lambda)$ on these measures, we then look for a $1$-homogeneous, convex, continuous function $\Psi_\alpha$ such that
\[
\Psi_\alpha(p) = ||g||_{1/\alpha} \quad \text{whenever } p \in K_\alpha = \{ \tau \otimes g, \quad ||\tau||_2=1, \; g_{i}\in\{0,1\} \}.
\]
The maximal function satisfying this condition can be computed as the $1$-homogeneous convex envelope of the function
\[
\Phi_\alpha(p) =
\left\{
\begin{aligned}
&||g||_{1/\alpha} &\quad& \text{if } p \in K_\alpha \\
&+\infty &\quad& \text{otherwise} \\
\end{aligned}
\right.
\]
and, as show in \cite{BoOrOu}, it turns out to be $\Phi_\alpha^{**}(p) = \sup_{q \in K^\alpha} \langle p, q \rangle$, which is to say the support function of the set
\[
K^\alpha=\left\{ p \in \R^{d\times(N-1)},\;\; \left\lVert \sum_{j\in J}p_j
\right\rVert_2 \le |J|^\alpha\ \ \forall\, J\subset\{1,...,N-1\}\ \right\},
\]
with $|J|$ the cardinality of the set $J$. Thanks to \eqref{eq:defpsilambda}, setting $\Psi_\alpha = \Phi_\alpha^{**}$, we can finally define
\[
\mathcal R^\alpha(\Lambda) = \sup \left\{\ \sum_{i=1}^{N-1} \int_{\R^d} \phi_i
\,d\Lambda_i,\quad \phi \in C_c^\infty (\R^d; K^\alpha) \right\},
\]
and consider the relaxed problem
\begin{equation}\label{eq:relaxed}
\inf \left\{  \mathcal{R}^\alpha(\Lambda), \quad \div \Lambda_i = \delta_{P_i}-\delta_{P_N} \text{ for all } i = 1,\dots,N-1 \right\}.
\end{equation}
This formulation provides the convex framework we were looking for: the problem is now defined on the whole
space of matrix valued Radon measures and
the energy is convex as it is a supremum of linear functionals.

However the functional $\mathcal R^\alpha$ is obtained only  as a ``local'' convex envelope of $\mathcal F^\alpha$ and as such it is not expected
 to always coincide with the true convex envelope, as we will see in Example \ref{ex:pentacol}. Thus, given a
 minimizer $\bar{\Lambda}$ of \eqref{eq:relaxed}
 we can end up in three different situations:
\begin{enumerate}
	\item $\bar{\Lambda} \in \mathcal{L}(A)$, then $\bar{\Lambda}$ is also a
     minimizer of $\mathcal{F}^\alpha$ and we have solved our original problem;
	\item $\mathcal{R}^\alpha(\bar{\Lambda}) = \inf_\Lambda \mathcal{F}^\alpha(\Lambda)$,
     then $\bar{\Lambda}$ is a convex combination of minimizers of $\mathcal{F}^\alpha$;
	\item $\mathcal{R}^\alpha(\bar{\Lambda}) < \inf_\Lambda \mathcal{F}^\alpha(\Lambda)$,
    which means that the relaxation is not tight and generally speaking minima of
    $\mathcal{R}^\alpha$ have no relation with minima of $\mathcal{F}^\alpha$.
\end{enumerate}
For a given set of terminal points $A = \{P_1,\dots,P_N\}$ we will then call the relaxation \eqref{eq:relaxed} to be tight (or sharp) whenever one of its minimizers satisfies 1. or 2., i.e. whenever its minimizers are related to the actual minimizers of $\mathcal{F}^\alpha$ as it is the case with real convex envelopes. Unfortunately, as the following counterexample shows, the relaxation is not
always sharp.

\begin{example}\label{ex:pentacol}[Non sharpness for pentagon configurations]
	Consider as terminal points the five vertices of a regular pentagon of side
	$\ell > 0$ and let $\beta = \frac{3}{10}\pi$. In this situation (STP) has $5$
	minimizers which are the one in the left picture of figure \ref{fig:steiner_vs_stella}
	and its $4$ rotations. The energy $\mathcal{R}^0$ of a Steiner tree, which
	corresponds by construction to its length, is equal to
	$\ell\tan\beta\, (1+\sin\beta+\sqrt{3}\cos\beta) \approx 3.8911\cdot \ell$.
	However none of the optimal Steiner trees is a minimizer for \eqref{eq:relaxed}.
	Indeed we can exhibit an admissible tensor valued measure $\Sigma$ with an
	energy strictly less than the energy of a Steiner tree: consider for example
	the rank one tensor valued measure $\Sigma$ constructed in the right picture
	of figure \ref{fig:steiner_vs_stella}. Such a measure satisfies the divergence
	constraints and its energy, which amounts to $1/2$ the length of its support,
	is equal to $\frac54\ell(\sqrt{3}+\tan\beta) \approx 3.8855 \cdot \ell$. Hence
	we are in the third case of the previous list: the relaxation is not tight and
	as we already said there is in general no way of reconstructing an optimum for
	(STP) staring from a minimizer of $\mathcal{R}^0$ (in this case our numerical
	results suggest $\Sigma$ as the actual minimizer of $\mathcal{R}^0$). Another
	example of non-sharpness can be obtained considering as terminal points the
	vertices of the pentagon plus the center: also in this case $\Sigma$ has less
	energy than any optimal Steiner tree.
	\begin{figure}
		\centering
		\begin{tabular}{cccc}
			\includegraphics[width=0.4\linewidth]{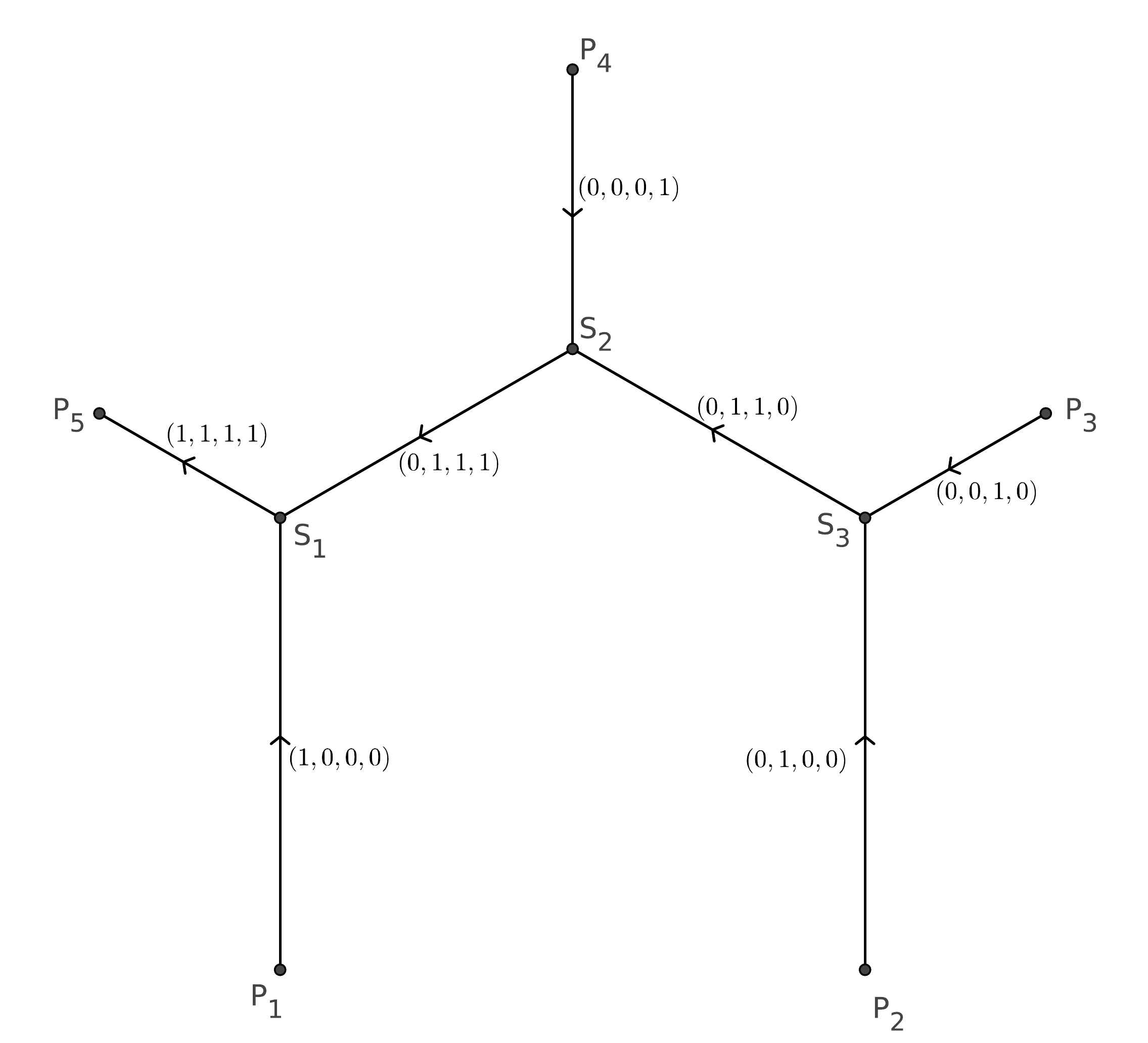} &
			\includegraphics[width=0.4\linewidth]{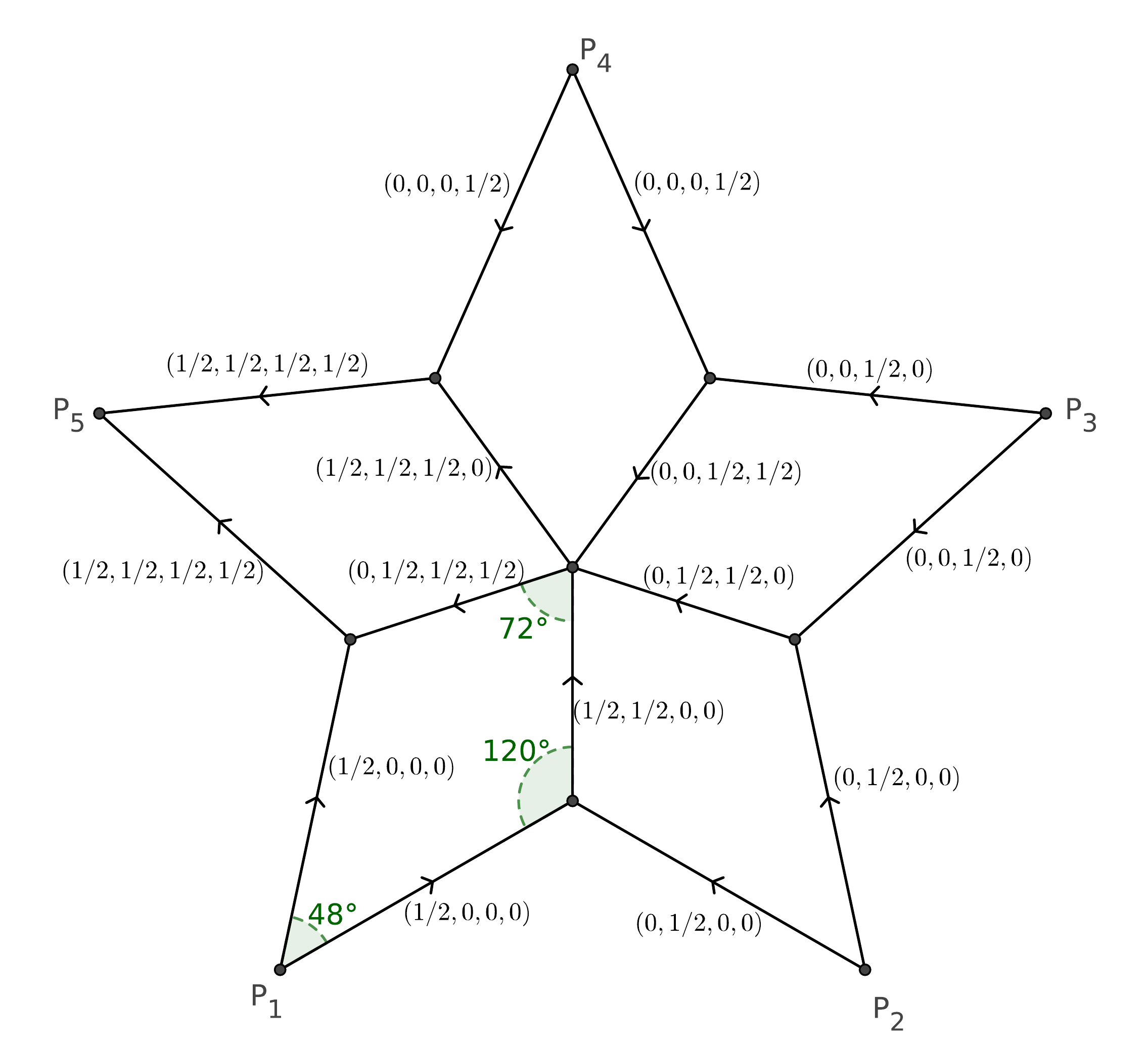}
		\end{tabular}
		\caption{Left: an optimal Steiner tree viewed as its corresponding measure $\Lambda$. Right: a rank one tenor valued measure $\Sigma=\tau \otimes g \cdot \mathcal{H}^1\res L$, with $L$ the graph itself, $\tau$ and $g$ as displayed.}
		\label{fig:steiner_vs_stella}
	\end{figure}
\end{example}

\smallskip


Despite the previous example, the proposed relaxation can be proved to be sharp in many situations. Indeed, thanks to the duality nature of $\mathcal{R}^\alpha$, we can prove minimality of certain given measures by means of calibration type arguments. This implies that whenever we are able to find a calibration for a given $\bar{\Lambda} \in \arg \min_\Lambda \mathcal{F}^\alpha(\Lambda)$ then the relaxation is sharp because $\bar{\Lambda}$ will also be a minimizer for $\mathcal{R}^\alpha$. A calibration, at least in the simple case of $\R^2$, can be defined as follows
\begin{definition}
	Fix a matrix valued Radon measure $\Lambda = (\Lambda_1, \dots, \Lambda_{N-1})$ and
	$\phi \in C^\infty_c(\R^2; K^\alpha)$. We say that $\phi$ is a calibration for
	$\Lambda$ if $\nabla \times \phi_i = 0$ for all $i = 1,\dots,N-1$, and
	$\phi$ realizes the supremum in the definition of $\mathcal R^\alpha$, i.e.
	\[
	\sum_{i=1}^{N-1} \int_{\R^2} \phi_i \,d\Lambda_i = \mathcal{R}^\alpha(\Lambda).
	\]
\end{definition}
The only existence of such an object certifies the optimality of $\Lambda$ in \eqref{eq:relaxed}. Indeed, let $\Sigma = (\Sigma_1, \dots,
\Sigma_{N-1})$ be another competitor, with $\mathcal{R}^\alpha(\Sigma) < \infty$ and $\div \Sigma_i = \delta_{P_i}-\delta_{P_N}$ for each $i = 1,\dots,N-1$. Hence $\div (\Lambda_i - \Sigma_i) = 0$ and we have\footnote{This generalizes the ``smooth'' case: thinking to $\Lambda_i$ and $\Sigma_i$ as ``regular'' vector fields we have that $\Lambda_i-\Sigma_i$ is a gradient, whence integrating by parts and using that $\phi_i$ is curl-free we get zero.}
\begin{equation}\label{eq:zeroint}
\int_{\R^2} \phi_i \,d(\Lambda_i - \Sigma_i) = 0,
\end{equation}
so that
\[
\begin{aligned}
\mathcal{R}^\alpha(\Lambda) &= \sum_{i=1}^{N-1} \int_{\R^2} \phi_i \,d\Lambda_
i = \sum_{i=1}^{N-1} \left( \int_{\R^2} \phi_i \,d(\Lambda_i-\Sigma_i) +
\int_{\R^2} \phi_i \,d\Sigma_i\right) \\
&\leq 0 + \mathcal{R}^\alpha(\Sigma) = \mathcal{R}^\alpha(\Sigma).
\end{aligned}
\]
In $\R^d$ with $d>2$,
the definition of a calibration extends as it is, where now
$\nabla \times \phi_i$ stands for the exterior derivative of the $1$-form
associated to the vector field $\phi_i$. Also \eqref{eq:zeroint} generalizes
and the proof carries over directly.

For the case $\alpha = 0$, which corresponds to (STP), we can take advantage of calibration arguments of \cite{MaMa} to justify sharpness of \eqref{eq:relaxed} for some classical choices of $\{P_1,\dots,P_N\}$. Indeed, as we will see in the next section, whenever $\Lambda$ is a
rank one tensor valued measure, for instance whenever it concentrates on a graph
and has real-valued weights, $\mathcal{R}^\alpha$ coincides
with the norm introduced in \cite{MaMa} to study (STP) as a mass-minimization problem for $1$-dimensional currents with
coefficients in a suitable normed group. Thus, every calibrated example in that context turns out to be a calibrated configuration in our framework, i.e. a situation where $\mathcal{R}^0$ is sharp (see \cite{MaMa, MaOuVe}).

\subsection{Extensions: generic Gilbert--Steiner problems and manifolds}

The same ideas developed in the previous paragraph can be extended beyond the (single sink) Gilbert--Steiner problem $(I_\alpha)$ in order
to address problems with possibly multiple sources/sinks in an Euclidean setting or even problems formulated within manifolds.

Following the strategy introduced in \cite{MaMa2} the energy
$\mathcal{F}^\alpha$ can also be used to address the general (oriented version of) ``who goes where''
problem. In this context we do not have to move all the mass to a single sink
but instead we are given a family of source/sink couples and we have to move a
unit mass from each source to each given destination. Thus, letting
$\{S_1,\dots,S_m\} \subset \R^d$ be the set of (unit) sources and
$\{T_1,\dots,T_m\} \subset \R^d$ the corresponding set of (unit) sinks, we
optimize the same energy $E^\alpha$ involved in the definition of $(I_\alpha)$
but this time among oriented networks of the form $L = \cup_{i=1}^m \lambda_i$, with
$\lambda_i$ a simple rectifiable curve connecting $S_i$ to $T_i$. The same
derivation as above can then be repeated, leading us to the relaxed formulation
\begin{equation}\label{eq:branched}
\inf \{ \mathcal{R}^\alpha(\Lambda), \quad \Lambda = (\Lambda_1,\dots,\Lambda_{m})\, , \div \Lambda_i = \delta_{S_i} - \delta_{T_i} \text{ for all } i = 1,\dots,m \}.
\end{equation}
We remark that in the previous who goes where problem, differently to what happens in \cite{BeCaMo}, we do not allow two paths $\lambda_i$, $\lambda_j$ to have opposite orientation on intersections, i.e. particles have to go the same way when flowing in the same region.

The previous approach to the ``who goes where'' problem can now be used within
the formulation of more general branched transportation problems, where we are
just required to move mass from a set of (unit) sources $\{S_1,\dots,S_m\}
\subset \R^d$ to a set of (unit) sinks $\{T_1,\dots,T_m\} \subset \R^d$,
without prescribing the final destination of each particle. In this context the
problem can be tackled as follows: for every possible coupling between sources
and sinks, i.e. among all permutations $\sigma \in \mathcal{S}_m$, solve the
corresponding ``who goes where'' problem with pairs $(S_i,
T_{\sigma(i)})_{i=1}^m$, and then take the coupling realizing the minimal
energy. Each ``who goes where'' can be relaxed as done in \eqref{eq:branched},
providing this way a relaxed formulation also for the case of generic multiple sources/sinks.

We point out how the extension of the previous discussion to a manifold framework is direct: the derivation that led us to the energy $\mathcal{R}^\alpha$, together
with problems \eqref{eq:relaxed} and \eqref{eq:branched}, is still valid on surfaces
embedded in the three dimensional space, with the only difference that divergence
constraints have to be intended as involving the tangential divergence operator
on the given surface.



\section{A first simple approximation on graphs}\label{sec:graph}

In this section we first see how the previous formulation simplifies when we consider the Steiner tree problem in the context of graphs, in which case the energy reduces to the norm introduced in \cite{MaMa}. Then,
once we are able to address (STP) on networks, we try to approximate the Euclidean (STP) by means of a discretization of the domain through an augmented graph.


\subsection{The Steiner tree problem on graphs}
Consider a connected graph $G = (V, E)$ in $\R^d$, where $V = \{v_1, \dots,
v_n\} \subset \R^d$ and $E = \{e_1, \dots,e_m\}$ is a set of $m$ segments. Each
$e_j = [v_{j}^1, v_{j}^2]$ connects two vertices $v_{j}^1, v_{j}^2$, has length
$\ell(e_j) = ||v_j^2-v_j^1||_2$ and is oriented by $\tau_{e_j} =
(v_{j}^2-v_{j}^1)/|v_{j}^2-v_{j}^1|$. Furthermore, we can assume without loss of
generality that edges intersect each other in at most $1$ point. The Steiner
Tree Problem within $G$ can be formulated in the same fashion as its Euclidean
counterpart: given a set of terminal points $A=\{P_1, \dots, P_N\} \subset V$
find the shortest connected sub-graph spanning $A$. As in the Euclidean case a
solution always exists and optimal sub-graphs are indeed sub-trees (they contain
no cycles).

Following what we did above in the Euclidean case, we can decompose any
candidate sub-graph $L \subset G$ into the superposition of $N-1$ paths
$\lambda_i$ within the graph, each one connecting $P_i$ to $P_N$. Each path is
identified as the support of a flow $V_i \colon E \to \{-1,0,1\}$ flowing a unit
mass from $P_i$ to $P_N$: we set $V_i(e)=1$ if edge $e$ is travelled in its own
direction within path $\lambda_i$, $-1$ if it is travelled in the opposite way and
$0$ otherwise. By construction we satisfy the discrete version of
\eqref{eq:solenoidal}, i.e. the classical Kirchhoff conditions: for all
``interior'' vertices $v \in V \setminus \{P_i, P_N\}$ we have
\begin{subequations}\label{eq:kirk}
\begin{equation}
\sum_{e\in \delta^+(v)} V_i(e) - \sum_{e\in \delta^-(v)} V_i(e) = 0,
\end{equation}
with $\delta^{\pm}(v)$ the set of outgoing/incoming edges at vertex $v$, and
$(P_i, P_N)$ is the source/sink couple, meaning
\begin{equation}
\sum_{e\in \delta^+(P_i)} V_i(e) - \sum_{e\in \delta^-(P_i)} V_i(e) = 1, \qquad
\sum_{e\in \delta^+(P_N)} V_i(e) - \sum_{e\in \delta^-(P_N)} V_i(e) = -1.
\end{equation}
\end{subequations}
Setting $V = (V_1,\dots,V_{N-1})$ and $L = \textup{supp}\,V = \cup
\{e \in E \,:\, V(e) \neq 0\}$, we have
\[
\mathcal{H}^1(L) = \sum_{e \in E} \ell(e) \cdot||V(e)||_\infty =: \mathcal{F}(V),
\]
and as before a solution to the network (STP) can be found minimizing
$\mathcal{F}$ among vector valued flows $V \colon E \to \{-1,0,1\}^{N-1}$
satisfying the above flux conditions \eqref{eq:kirk}. Let us identify each
family $V$ with a tensor valued measure $\Lambda = (\Lambda_1, \dots,\Lambda_{N-1})$
defined on the whole $\R^d$ by setting
\begin{equation}\label{eq:fromVtoL}
\Lambda_i = \sum_{e \in E} V_i(e)\, \tau_e\cdot \mathcal{H}^1\res e, \quad
\text{i.e.} \quad \Lambda = \sum_{e \in E} \tau_e \otimes V(e) \cdot \mathcal{H}^1\res e.
\end{equation}
The idea is now to drop the integer constraint $\{-1,0,1\}$ on each $V_i$ and
optimize the previously defined energy $\mathcal{R}^0$ among tensor valued measures of the form
\eqref{eq:fromVtoL}, obtaining the relaxed energy
\[
\mathcal{R}(V) = \mathcal{R}^0(\Lambda) = \sup_{\phi \in C_c^\infty (\R^d; K^0)}
\ \sum_{i=1}^{N-1} \int_{\R^d} \phi_i \,d\Lambda_i = \sup_{\phi \in C_c^\infty
(\R^d; K^0)} \ \sum_{i=1}^{N-1} \sum_{e \in E} \left(V_i(e) \int_{e} \phi_i \,ds\right).
\]
Since edges intersect in at most $1$ point it is possible to interpret the last
supremum as a supremum over test functions entirely supported on the graph and
of the form $\phi = \sum_{e} \tau_e \otimes W(e)$ with $W \colon E \to \R^{N-1}$.
By assumption, for almost every point $x$ on the graph (except at intersections)
there exists only one edge $e$ containing $x$; hence, the pointwise constraint
$\phi(x) \in K^0$ translates into $\phi\res e \in K^0$ for all edges $e \in E$, i.e.
\[
\left\lVert\sum_{j\in J}W_j(e)\tau_e \right\rVert_2 = \left|\sum_{j\in J}W_j(e)
\right|\le 1\quad \forall\, J\subset\{1,...,N-1\}.
\]
These new constraints involve only vectors $W(e)$ and are equivalent to the
unique constraint
\[
||W(e)||_* = \left[\sum_{j=1}^{N-1} (W_j(e) \vee 0)\right] \vee \left[ -
\sum_{j=1}^{N-1} (W_j(e) \wedge 0) \right] \leq 1,
\]
which amounts to require that the maximum between the $\ell^1$ norm of the positive
part and the $\ell^1$ norm of the negative part of $W(e)$ has to be less or
equal to $1$. The energy can be finally rewritten as
\[
\begin{aligned}
\mathcal{R}(V) &= \sup \left\{ \sum_{e\in E} \ell(e) \,V(e) \cdot W(e), \quad
||W(e)||_* \leq 1 \;\forall e \in E \right\}\\
&= \sum_{e \in E} \ell(e)\left( \sup_i [V_i(e) \vee 0] - \inf_i [V_i(e) \wedge 0]
\right) = \sum_{e \in E} \ell(e)||V(e)||.
\end{aligned}
\]
The norm $||\cdot||$ is exactly the norm used in \cite{MaMa} to study (STP) using
currents with coefficients in normed groups and hence we can take advantage
of calibration arguments of \cite{MaMa} to justify the sharpness of the relaxation for calibrated configurations of terminal points. Of course the counterexample \ref{ex:pentacol} still applies to this discrete version of the problem
using as graph $G$ the union of the two graphs of picture \ref{fig:steiner_vs_stella}:
the minimizer concentrates on the star and not on the Steiner
structure.

Optimization of $\mathcal{R}$ under the (linear) flux constraints \eqref{eq:kirk} can then be performed solving a linear program: in order to linearize the objective we introduce two sets of variables $\{s_e\}_{e\in E}$, $\{i_e\}_{e\in E}$, and for each $e \in E$ we require $i_e \leq 0$, $s_e \geq 0$ and
\[
i_e \leq V_i(e) \leq s_e \quad \text{for all }i=1,\dots,N-1,
\]
so that the objective reduces to $\sum_e \ell(e) (s_e-i_e)$. Whenever the size of the resulting linear program is too big to be treated by standard interior point solvers we can alternatively apply the cheaper first order scheme proposed in \cite{ChPo} (see Section \ref{sec:algo} for details).


\subsection{Graphs and the Euclidean (STP)}

\begin{figure}
	\centering
	\begin{tabular}{ccc}
		\includegraphics[width=0.29\linewidth]{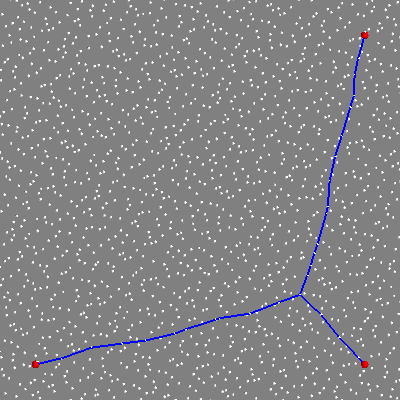} &
		\quad
		\includegraphics[width=0.29\linewidth]{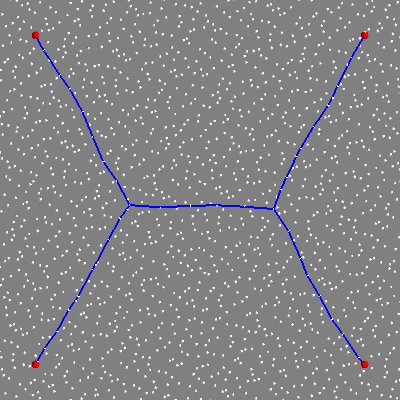} &
		\quad \includegraphics[width=0.29\linewidth]{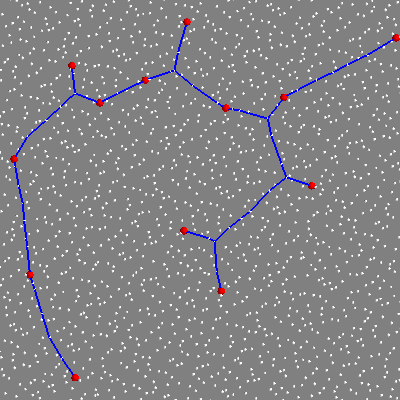}
	\end{tabular}
	\caption{Approximations of (STP) for $3$, $4$ and $13$ terminal points (red) using the augmented graph idea, $K = 1681$, $M = 30$. Edges carrying a non-zero flux are displayed.}
	\label{fig:graph_steiner}
\end{figure}

Once we have a method to approximate (STP) on networks we can try to address the Euclidean (STP) through the use of an augmented graph. The core idea is the following: let $\{x_1, \dots, x_K\}$, $K \in \mathbb{N}$, be a set of scattered points that uniformly covers an open convex domain $\Omega$ such that $\{P_1,\dots,P_N\} \subset \Omega$ and let $V = \{ x_1, \dots, x_K\} \cup \{ P_1, \dots, P_N \}$. Fix $M \in \mathbb{N}$ and construct the graph $G = (V, E)$ where each $v \in V$ is connected through segments to its $M$ closest neighbours. For $M$ sufficiently large the network $G$ is connected and solving (STP) within $G$ provides an approximation of the underlying Euclidean Steiner tree.

We see in figure \ref{fig:graph_steiner} two examples with $K = 1681$ and $M = 30$. In both cases results are very close to the optimal Steiner tree and for obtaining them we simply solve a medium scale linear program. However the use of a fixed underlying graph has some drawbacks. For example we cannot expect edges meeting at triple points to satisfy the $120^\circ$ condition and what should be a straight piece in the optimal tree is only approximated by a sequence of (non-aligned) edges. A possible remedy for obtaining ``straighter'' solutions is to increase $M$, allowing this way longer edges, but this would increase the size of the problem. Furthermore obtaining convex combinations of minimizers is almost impossible because the underlying graph is not regular and having two sub-graphs with the exact same energy is very rare. On the other hand taking regularly distributed points generates many equivalent solutions even when there should be only one.

We also observe that this simplified framework is specific to the Euclidean Steiner tree case: the corresponding graph framework for $(I_\alpha)$ does not end up in a linear program and no direct extension to the manifold case is possible. This lack of generality, together with the intrinsic low precision of the approach as a consequence of working on a graph, leads us to switch our focus on the direct minimization of $\mathcal{R}^\alpha$ on the whole of $\R^2$/$\R^3$.

\section{Generic Euclidean setting, the algorithmic approach}\label{sec:algo}
Motivated by the shortcomings of the previous simplified framework, we present in this section our approach for solving \eqref{eq:relaxed} in $\R^2$ (the same ideas extends to the three dimensional setting). Our resolution is based on a staggered grid for the discretization of the unknowns coupled with a conic solver (or a primal-dual scheme) for the optimization of the resulting finite dimensional problem.

\subsection{Spatial discretization}
Assume without loss of generality that $P_1,\dots,P_N$ are contained in the interior of $\Omega = [0,1]\times [0,1]$, which will be our computational domain. From a discrete standpoint we view the unknown vector measures $(\Lambda_1, \dots, \Lambda_{N-1})$ in \eqref{eq:relaxed} as a family $V = (V_1, \dots,
V_{N-1})$ of vector fields in $\Omega$ and, due to the divergence constraints
that we need to satisfy, we discretize these unknown fields on a staggered grid
(this way our degrees of freedom are directly related to the flux of each vector
field through the given grid interface). Fix then a regular Cartesian grid of
size $M\times M$ over $\Omega$ and let $h = 1/M$.
The first component $V_{i,1}$ of each vector field is placed on the midpoints of
the vertical cells interfaces whereas the second components $V_{i,2}$ on the
horizontal ones, so that to have on each element $(k,\ell)$ \[
V_{i}|_{(k,\ell)} = \begin{pmatrix}
(V_{i,1}^{k+1,\ell}-V_{i,1}^{k,\ell}) (x-(k-1)h)/h + V_{i,1}^{k,\ell} \\ (V_{i,2}^{k,\ell+1}-V_{i,2}^{k,\ell}) (y-(\ell-1)h)/h + V_{i,2}^{k,\ell}
\end{pmatrix}.
\]
The component $V_{i,1}$ is described by $(M + 1) \times M$ unknowns whereas
$V_{i,2}$ is described by $M \times (M + 1)$ parameters. Regarding the test
functions $\phi = (\phi_1, \dots, \phi_{N-1})$ we define them to be piecewise
constant on each element of the grid, i.e. for any cell $(k,\ell)$ we have
$\phi_i^{k,\ell} = (\phi_{i,1}^{k,\ell}, \phi_{i,2}^{k,\ell}) \in \R^2$.

Within this setting the optimization of the energy $\mathcal{R}^\alpha$ translates
into
\begin{equation}\label{eq:discrete_energy}
\min_{(V_{i,d}^{k,\ell})}\sup_{(\varphi_{i,d}^{k,\ell}) \in K^\alpha} \sum_{k,\ell} \sum_{i=1}^{N-1} h^2 \, \left[\frac{V_{i,1}^{k,\ell} + V_{i,1}^{k+1,\ell}}{2}\,\varphi_{i,1}^{k,\ell} + \frac{V_{i,2}^{k,\ell} + V_{i,2}^{k,\ell+1}}{2}\, \varphi_{i,2}^{k,\ell}\right]
\end{equation}
under the condition $\div V_i = \delta_{P_i}-\delta_{P_N}$ for all
$i=1,\dots,N-1$. Since the flux of each $V_i$ over the generic cell $(k,\ell)$
is given by
\[
F_i^{k,\ell} = h(V_{i,1}^{k+1,\ell} - V_{i,1}^{k,\ell}) + h(V_{i,2}^{k,\ell+1} - V_{i,2}^{k,\ell}),
\]
the divergence constraints translate, at a discrete level, into
\begin{equation}\label{eq:discretediv}
\left\lbrace
\begin{aligned}
F_i^{k,\ell} &= 0 &\quad& \text{whenever cell $(k,\ell)$ does not contain $P_i$ or $P_N$}, \\
F_i^{k,\ell} &= 1 &\quad& \text{if cell $(k,\ell)$ contains $P_i$}, \\ F_i^{k,\ell}
&= -1 &\quad& \text{if cell $(k,\ell)$ contains $P_N$},
\end{aligned}
\right.
\end{equation}
complemented with a ``zero flux'' condition at the boundary, i.e. we set
$V_{i,d}^{k,\ell} = 0$ whenever it refers to a boundary interface. We finally
observe that, by construction, $\phi \in K^\alpha$ if for each cell $(k,\ell)$
in the grid the matrix $\phi^{k,\ell}
=(\phi^{k,\ell}_1,\dots,\phi^{k,\ell}_{N-1})$ satisfies
\[
\left\lVert \sum_{j\in J} \varphi_j^{k,\ell} \right\lVert_2 \leq |J|^\alpha \text{ for all }J\subset \{ 1,\dots,N-1 \}.
\]
For the resolution of this finite dimensional optimization problem we then
propose two different and somehow complementary approaches.
\subsection{Optimization via conic duality}
The $\inf$-$\sup$ problem \eqref{eq:discrete_energy} can be written, thanks to
conic duality (see e.g. Lecture 2 of \cite{ben2001lectures}), as a pure
minimization problem involving the degrees of freedom $(V_{i,d}^{k,\ell})$ and
a set of dual variables $(\psi_{J,d}^{k,\ell})$ indexed over subsets
$J\subset\{1,\dots,N-1\}$. Indeed, for fixed $1\leq k,\ell\leq M$ and
 $J \subset \{1,\dots,N-1\}$, one has
\[
\inf_{\psi_J^{k,\ell} \in \R^2} \left(|J|^\alpha||\psi_J^{k,\ell}||_2 - \langle \,
\psi_J^{k,\ell}, \sum_{j\in J} \varphi_j^{k,\ell} \,\rangle\right) =
\left\lbrace
\begin{aligned}
&0 &\quad& \text{ if } \left\lVert \sum_{j\in J} \varphi_j^{k,\ell} \right\lVert_2 \leq |J|^\alpha \\
&-\infty &\quad& \text{ otherwise,} \\
\end{aligned}
\right.
\]
so that, if we denote $\tilde{V}_i^{k,\ell} = ((V_{i,1}^{k,\ell} +
V_{i,1}^{k+1,\ell})/2, (V_{i,2}^{k,\ell} + V_{i,2}^{k,\ell+1})/2) \in \R^2$,
\eqref{eq:discrete_energy} is equivalent to
\[
\min_{(V_{i,d}^{k,\ell})}\sup_{(\varphi_{i,d}^{k,\ell})} \left[ \sum_{k,\ell}
\sum_{i=1}^{N-1} h^2 \langle \tilde{V}_i^{k,\ell}, \phi_i^{k,\ell}\rangle +
\inf_{(\psi_{J,d}^{k,\ell})} \sum_{k,\ell}\sum_{J} h^2\left(|J|^\alpha||\psi_J^{k,\ell}||_2
- \langle \,\psi_J^{k,\ell}, \sum_{j\in J} \varphi_j^{k,\ell} \,\rangle\right)\right].
\]
Switching the $\sup$ over $(\varphi_{i,d}^{k,\ell})$ and the $\inf$ over $(\psi_{J,d}^{k,\ell})$ we obtain
\[
\min_{(V_{i,d}^{k,\ell})}\inf_{(\psi_{J,d}^{k,\ell})} \left[ \sum_{k,\ell}
\sum_{J} h^2|J|^\alpha||\psi_J^{k,\ell}||_2 + h^2\sup_{(\varphi_{i,d}^{k,\ell})}
\sum_{k,\ell}\left(\sum_{i=1}^{N-1} \langle \tilde{V}_i^{k,\ell}, \phi_i^{k,\ell}\rangle -
\langle \,\psi_J^{k,\ell}, \sum_{j\in J} \varphi_j^{k,\ell} \,\rangle\right)\right].
\]
Since the inner $\sup$ is either $0$ if $\tilde{V}_i^{k,\ell} = \sum_{J\ni i} \psi_J^{k,\ell}$
for all $1\leq k,\ell\leq M$ and $1\leq i\leq N-1$ or $+\infty$ otherwise,
the previous problem eventually leads to
\begin{equation}\label{eq:psienergy}
\min_{(V_{i,d}^{k,\ell}), (\psi_{J,d}^{k,\ell})} \sum_{k,\ell} \sum_{J} h^2 |J|^\alpha\,\norm{\psi_{J}^{k,\ell}}_2
\end{equation}
where each $V_i$ satisfies the same flux constraints \eqref{eq:discretediv} and
 for all cells $(k,\ell)$ and all $i=1,\dots,N-1$ we must satisfy
\begin{equation}\label{eq:Vpsirelation}
\frac{V_{i,1}^{k,\ell} + V_{i,1}^{k+1,\ell}}{2} = \sum_{J \ni i} \psi_{J,1}^{k,\ell}
\quad \text{and} \quad \frac{V_{i,2}^{k,\ell} + V_{i,2}^{k,\ell+1}}{2} = \sum_{J \ni i} \psi_{J,2}^{k,\ell}.
\end{equation}
Problem \eqref{eq:psienergy} under the set of linear constraints \eqref{eq:discretediv}
and \eqref{eq:Vpsirelation} can now be solved invoking the conic solver of the library
MOSEK \cite{mosek2010mosek} within the framework provided by \cite{JuMP}.

\subsection{Optimization via primal-dual schemes}


Collect all the $(V_{i,d}^{k,\ell})$ into a vector $\mathbf{v} \in \R^{n_v}$, $n_v = (N-1)(2M^2+2M)$, and all the $(\varphi_{i,d}^{k,\ell})$ into $\bm{\varphi} \in \R^{n_\phi}$, $n_\phi = (N-1)2M^2$. Moving the constraints on $\phi$ into the objective via the convex indicator function, the discrete energy \eqref{eq:discrete_energy} can be written down as
\[
\min_{\mathbf{v}} \max_{\bm{\phi}} \, \langle \bm{\phi}\,, B\mathbf{v} \rangle - \chi_{K^\alpha}(\bm{\phi})
\]
for a suitable (sparse) matrix $B$ of size $n_\phi \times n_v$, while the divergence constraints reduce to $A\mathbf{v} = \mathbf{b}$ for a suitable (sparse) matrix $A$ of size $n_\lambda \times n_v$ and a vector $\mathbf{b}\in \R^{n_\lambda}$. To the set of liner constraints $A\mathbf{v} - \mathbf{b} = 0$ we can now associate a dual variable $\bm{\lambda} \in \R^{n_\lambda}$ so that they can be incorporated into the objective as
\[
\min_{\mathbf{v}} \max_{\bm{\phi},\bm{\lambda}} \, \langle \bm{\phi}\,, B\mathbf{v} \rangle - \chi_{K^\alpha}(\bm{\phi}) + \langle \bm{\lambda}\,,A\mathbf{v} - b\rangle.
\]

The problem, written this way, turns into an instance of a general $\inf$-$\sup$ problem of the form
\begin{equation}\label{eq:generalpd}
\min_{x\in\R^n} \max_{y\in\R^m} \, \langle y\,, Kx \rangle + G(x) - F^*(y)
\end{equation}
with $K$ an $m \times n$ matrix and $G\colon \R^n \to \R \cup \{\infty\}, F^*\colon \R^m \to \R \cup \{\infty\}$ convex lsc functions. Among the possible numerical schemes which have been developed in the literature for the resolution of \eqref{eq:generalpd} we choose here the preconditioned primal-dual scheme presented in \cite{ChPo}. The scheme can be summarized as follows: let $\gamma\in[0,2]$, $T = \textup{diag}(\tau_1,\dots,\tau_n)$ and $\Sigma = \textup{diag}(\sigma_1,\dots,\sigma_m)$, with
\[
\tau_j = \frac{1}{\sum_{i=1}^{m}|K_{ij}|^{2-\gamma}} \quad \text{and} \quad \sigma_i = \frac{1}{\sum_{j=1}^{n}|K_{ij}|^{\gamma}},
\]
fix $x^0 \in \R^n$, $y^0 \in \R^m$, and iterate for any $k>0$
\begin{equation}\label{eq:pdscheme}
\left\lbrace
\begin{aligned}
x^{k+1} &= (I+T\partial G)^{-1}(x^k - TK^Ty^k) \\
y^{k+1} &= (I+\Sigma\partial F^*)^{-1}(y^k + \Sigma K(2x^{k+1}-x^k)) \\
\end{aligned}
\right.
\end{equation}
In this context the proximal mappings are defined as
\[
(I+T\partial G)^{-1}(\hat{x}) = \arg \min_x \left[G(x) + \frac12 \langle T^{-1}(x-\hat{x}), \,x-\hat{x} \rangle\right]
\]
and represent the extension of the classical definition with constant step size to this situation with ``variable dependent'' step sizes.

In our specific use case the scheme takes the following form: define $T = \textup{diag}(\tau_1,\dots,\tau_{n_v})$, $\Sigma = \textup{diag}(\sigma_1,\dots,\sigma_{n_\phi})$ and $\tilde{\Sigma} = \textup{diag}(\tilde{\sigma}_1,\dots,\tilde{\sigma}_{n_\lambda})$, with
\[
\tau_j = \frac{1}{\sum_{i=1}^{n_\phi}|B_{ij}|^{2-\gamma}+\sum_{i=1}^{n_\lambda}|A_{ij}|^{2-\gamma}}, \quad \sigma_i = \frac{1}{\sum_{j=1}^{n_v}|B_{ij}|^{\gamma}}, \quad \tilde{\sigma}_i = \frac{1}{\sum_{j=1}^{n_v}|A_{ij}|^{\gamma}},
\]
given $\mathbf{v}^{0}, \bm{\phi}^{0}, \bm{\lambda}^0$ iterate for $k>0$
\begin{equation}\label{eq:pd}
\left\lbrace
\begin{aligned}
\mathbf{v}^{k+1} &= \mathbf{v}^k - T(B^T\bm{\phi}^k + A^T\bm{\lambda}^k) \\
\bm{\phi}^{k+1} &= \textup{proj}(\bm{\phi}^{k} + \Sigma B(2\mathbf{v}^{k+1}-\mathbf{v}^{k}) \;|\; K^\alpha) \\
\bm{\lambda}^{k+1} &= \bm{\lambda}^k + \tilde{\Sigma}(A(2\mathbf{v}^{k+1}-\mathbf{v}^{k}) - \mathbf{b})
\end{aligned}
\right.
\end{equation}
The computational bottleneck for this simple iterative procedure resides in the
projection of a given vector $\bar{\bm{\phi}} \in \R^{n_\phi}$ onto the convex
set $K^\alpha$. By definition this operation reduces to the cell-wise projection
on $K^\alpha$ of the matrices $\phi^{k,\ell}$, and so we fix a $d \times (N-1)$
matrix $q = (q_1, \dots, q_{N-1})$ and split the discussion into two sub-steps.

\textbf{Projection on individual sets}: for each fixed subset $J \subset \{1,\dots,N-1\}$ we define the convex set
\[
K^\alpha_J = \left\{ p \in \R^{d\times(N-1)},\;\; \left\lVert\sum_{j\in J}p_j \right\rVert_2\le |J|^\alpha \right\}.
\]
The projection of $q$ over $K^\alpha_J$ can be computed explicitly: define $v_J =
\sum_{j \in J} q_j$, then the projection $p = \textup{proj}(q \;|\; K^\alpha_J) =
(p_1,\dots, p_{N-1})$ has columns defined as $p_j = q_j$ if $j \notin J$ and \[
p_j = q_j - {1}/{|J|}\,(\norm{v}_2-|J|^\alpha)^+ \frac{v}{\norm{v}_2} \quad \text{if $j \in J$.}
\]

\textbf{Projection on the intersection}: observe that $K^\alpha = \cap_J K^\alpha_J$, i.e. $K^\alpha$ is the intersection of a family of convex sets. In order to get an approximation of $\textup{proj}(q \;|\; K^\alpha)$ we can apply the Dykstra's projection algorithm (see \cite{dykstra83}). The scheme in our setting is the following: let $J_1, \dots, J_{2^{N-1}}$ be all the subsets of $\{1,\dots,N-1\}$, let $\{y_j^0\}_{j=1}^{2^{N-1}}$ be $2^{N-1}$ null matrices of size $d \times (N-1)$, $p^0 = q$, then for any $k\geq 1$ iterate
\[
\left\lbrace
\begin{aligned}
&p_0^{k} = p^{k-1} \\
&\text{for } j = 1,\dots,2^{N-1} \\
&\qquad\quad
\begin{aligned}
p_j^{k} &= \textup{proj}(p_{j-1}^{k} + y_{j}^{k-1} \;|\; K_{J_j}^\alpha) \\
y_j^{k} &= y_j^{k-1} + p_{j-1}^{k} - p_j^{k} \\
\end{aligned} \\
&\text{end for} \\
&p^{k} = p_{2^{N-1}}^{k}
\end{aligned}
\right.
\]
We then have $p^k \to \textup{proj}(q \;|\; K^\alpha)$ as $k \to +\infty$.

\begin{remark}
Each step of the previous iterative projection procedure requires $2^{N-1}$
sub-projections and thus the scheme is intrinsically time-consuming. Up to our knowledge
there seems to be no immediate simplifications to avoid some of the $2^{N-1}$
inner projections: for example the restriction of the inner loop over sets
$K^\alpha_{J_j}$ such that $q \notin K^\alpha_{J_j}$ is not going to work in
general. At the same time we observe that established convergence rates for
\eqref{eq:pdscheme} do not apply in this case because our projection, which
represents one of the two proximal mappings, is only approximated and not exact,
making us falling back in a context like \cite{ScLeBa}.

\end{remark}


\section{Numerical details}\label{sec:refinement}
The two resolution paths presented above allow us to overcome some shortcomings of the simplified framework of Section \ref{sec:graph} but introduces at the same time an higher computational cost, mainly depending on the combinatorial nature of the set $K^\alpha$, which reflects in the high number of variables involved in \eqref{eq:psienergy} and in the complicated projection in \eqref{eq:pd}.

\begin{figure}[t]
	\centering
	\begin{tabular}{cccc}
		\includegraphics[width=0.36\linewidth]{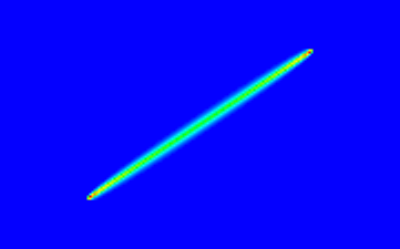} &
		\includegraphics[width=0.36\linewidth]{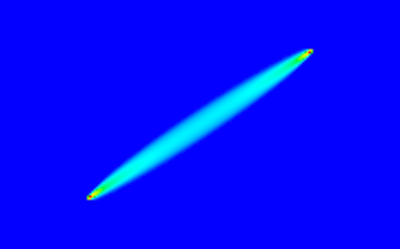}
	\end{tabular}
	\caption{Energy concentration for the minimizer of $\mathcal{R}^0$ for $P_1 = (1/4,1/3)$, $P_2 = (3/4, 2/3)$. Left: solution obtained via a conic solver, final energy $\approx 0.606307$. Right: solution obtained after $200000$ iterations of the primal-dual scheme, $\gamma = 0.6$, final energy $\approx 0.606765$.}
	\label{fig:conic_vs_primaldual}
\end{figure}

Generally speaking the primal-dual scheme is the cheapest of the two in terms of computational resources: it can be implemented so that every operation is done in-place, reducing to almost zero any further memory requirement apart from initialization, while the interior point approach used by a conic solver is extremely demanding in terms of memory due to the $2^{N-1}$ additional variables needed to define \eqref{eq:psienergy}. However, since we are looking for $1$d structures, our solver also needs to be able to provide very localized optima, and with this regards the primal-dual approach is not very satisfactory. As we can see in figure \ref{fig:conic_vs_primaldual}, where we use the two schemes for the same regular $201\times201$ grid over $[0,1]^2$, the solution provided by the primal-dual scheme is more diffused than the one obtained using the conic approach. For this reason we would like to use for our experiments the conic formulation \eqref{eq:psienergy} and to do so, in order to be able to treat medium scale problems, we need to find a way to reduce a-priori the huge number of additional variables that are introduced: this can be done both via a classical grid refinement and via a variables ``selection''.




\subsection{Grid refinement}
\begin{figure}[tbh]
	\centering
	\begin{tabular}{cccc}
		\includegraphics[width=0.219\textwidth]{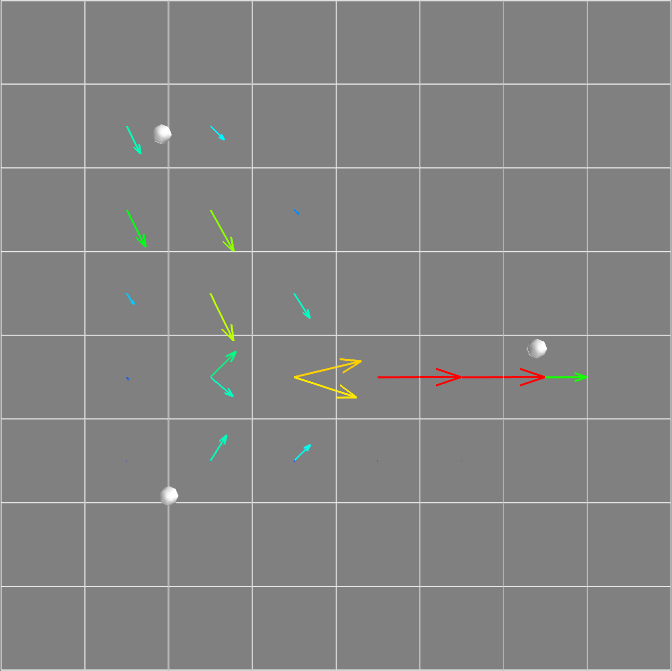} &
		\includegraphics[width=0.219\textwidth]{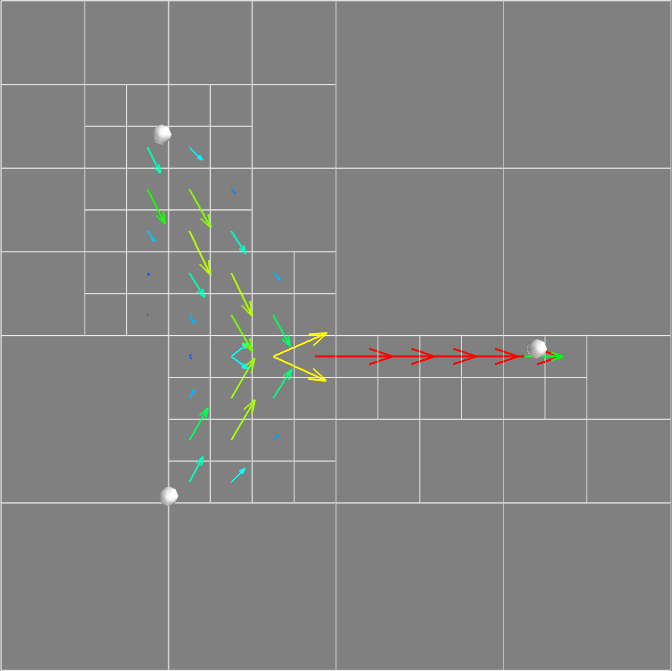} &
		\includegraphics[width=0.219\textwidth]{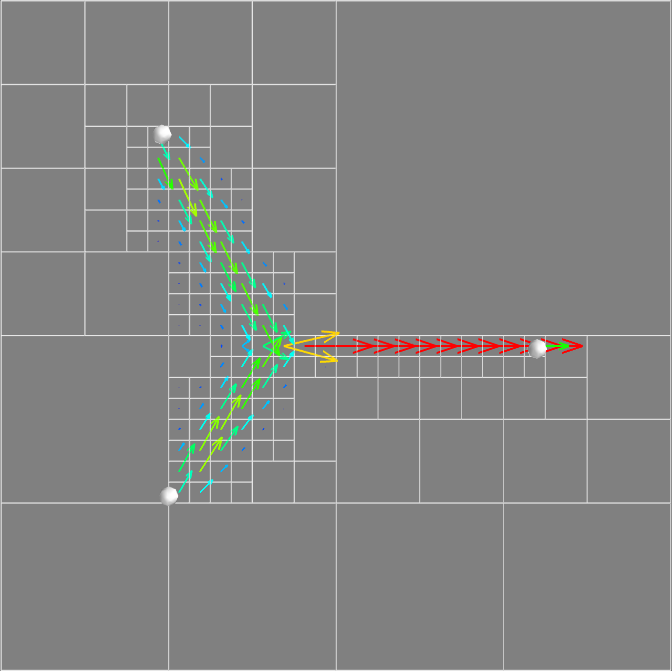} &
		\includegraphics[width=0.219\textwidth]{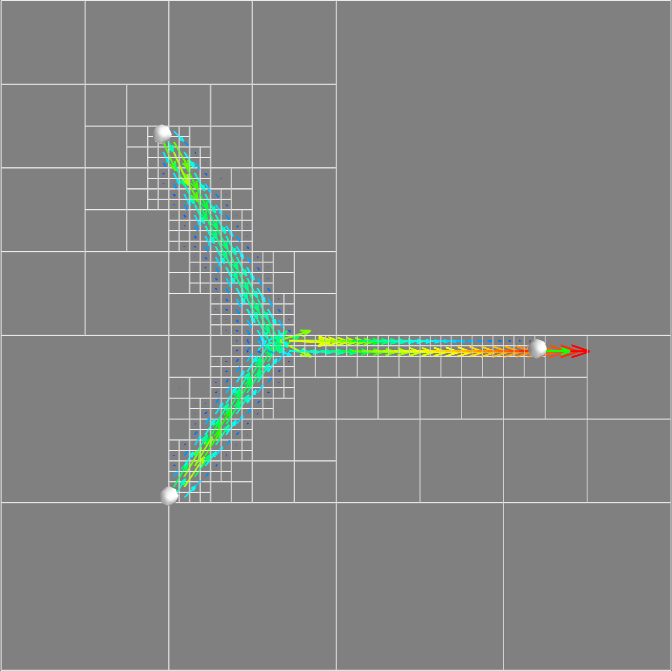}
	\end{tabular}
	\caption{Refinement example for $3$ points. At each iterate we plot the grid and the two fields $V_1$, $V_2$, which are then used to build the next grid.} 
	\label{fig:refinement}
\end{figure}

The numerical solution is expected to concentrate on a $1$-dimensional
structure, and so the grid needs to be fine only on a relatively small region of
the domain. This suggests the implementation of a refinement strategy able to
localize in an automatic way the region of interest. For doing so we use
non-conformal quadtree type meshes (see e.g. \cite{samet88, bass87}), which are
a particular class of grids where the domain is partitioned using $M$ square
cells as $\Omega = \cup_{m} S_m$ and each square cell $S_m$ can be obtained by
recursive subdivision of the box $[0,1]^2$ (see figure \ref{fig:refinement} for
examples of such grids). As in the case of uniform regular meshes we employ for
the discretization a staggered approach: we set the degrees of freedom of vector
fields on faces of each element, with the additional requirement that whenever a
face is also a subsegment of another face then the two associated degrees of
freedom are equal (this is to maintain continuity of the normal components of
the discrete fields across edges and guarantees that fluxes are globally well
behaved). The matrix valued function $\phi$ is again defined to be constant on
each element of the grid so that the nature of the discrete problem we need to
solve remains the same. 

A refinement procedure can then be described as follows: fix a coarse quadtree
grid $\mathcal{T}$, for example a regular $8\times 8$ one, and then
\begin{itemize}
	\item solve the problem on the given grid $\mathcal{T}$;
	\item identify elements of the grid where the solution concentrates the most
    and label them as ``used'', identify elements of the grid where the solution
    is almost zero and label them as ``unused'';
	\item refine the grid subdividing each ``used'' element into $4$ equal
    sub-elements and try to merge ``unused'' elements into bigger ones (the
    merging will occur if four elements labelled as ``unused'' have the same
    father in the quadtree structure);
	\item repeat.
\end{itemize}
As we can see in figure \ref{fig:refinement} this procedure allows us to localize
computations in a neighbourhood of the optimal structure we are looking for. This
way we can attain a good level of fineness around the solution without being forced
to employ a full grid which would require the introduction of a lot of useless degrees of freedom.

\subsection{Variables selection}

\begin{figure}[tb]
	\centering
	\begin{tabular}{cc}
		\includegraphics[width=0.45\textwidth]{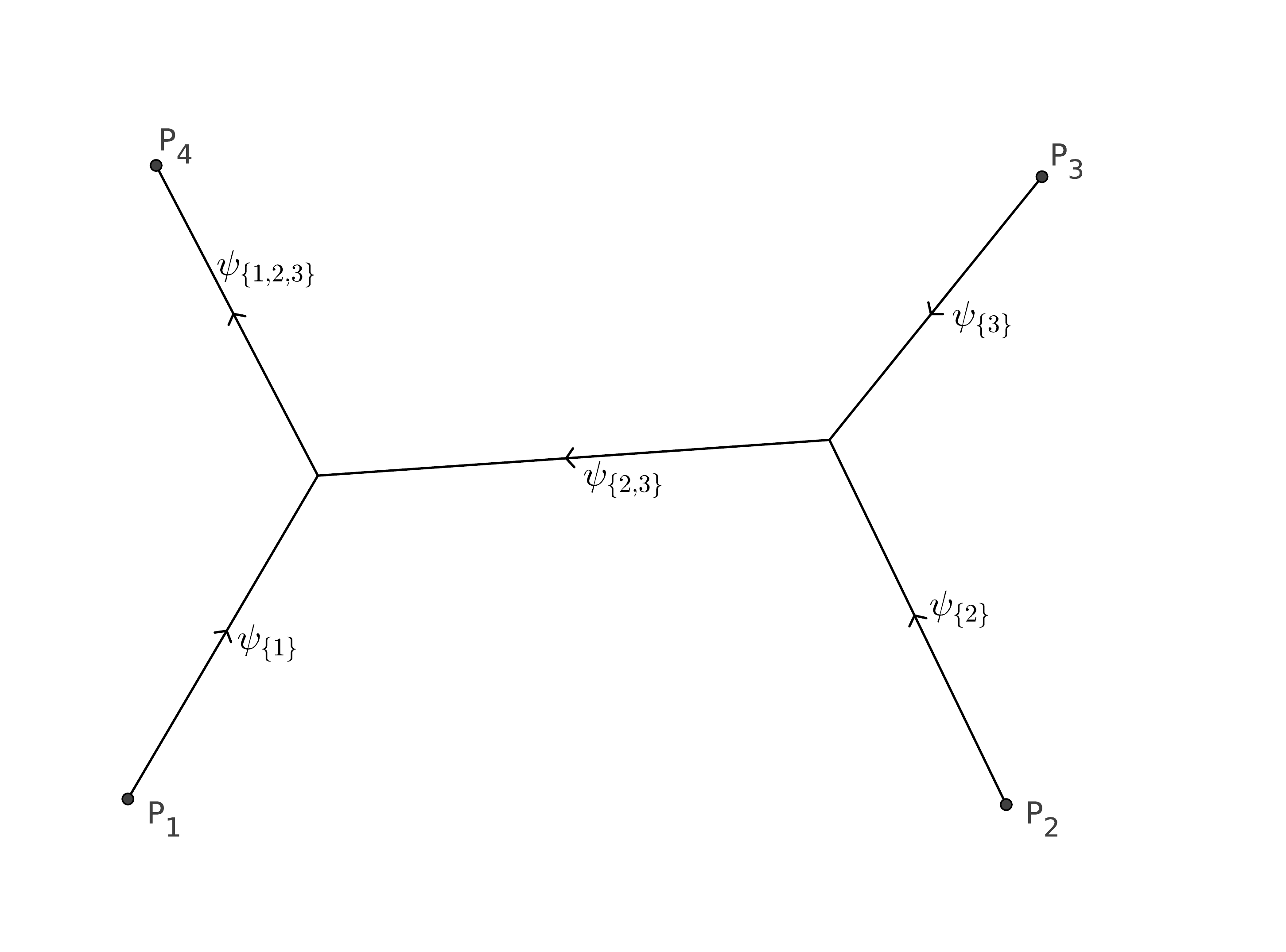} &
		\includegraphics[width=0.45\textwidth]{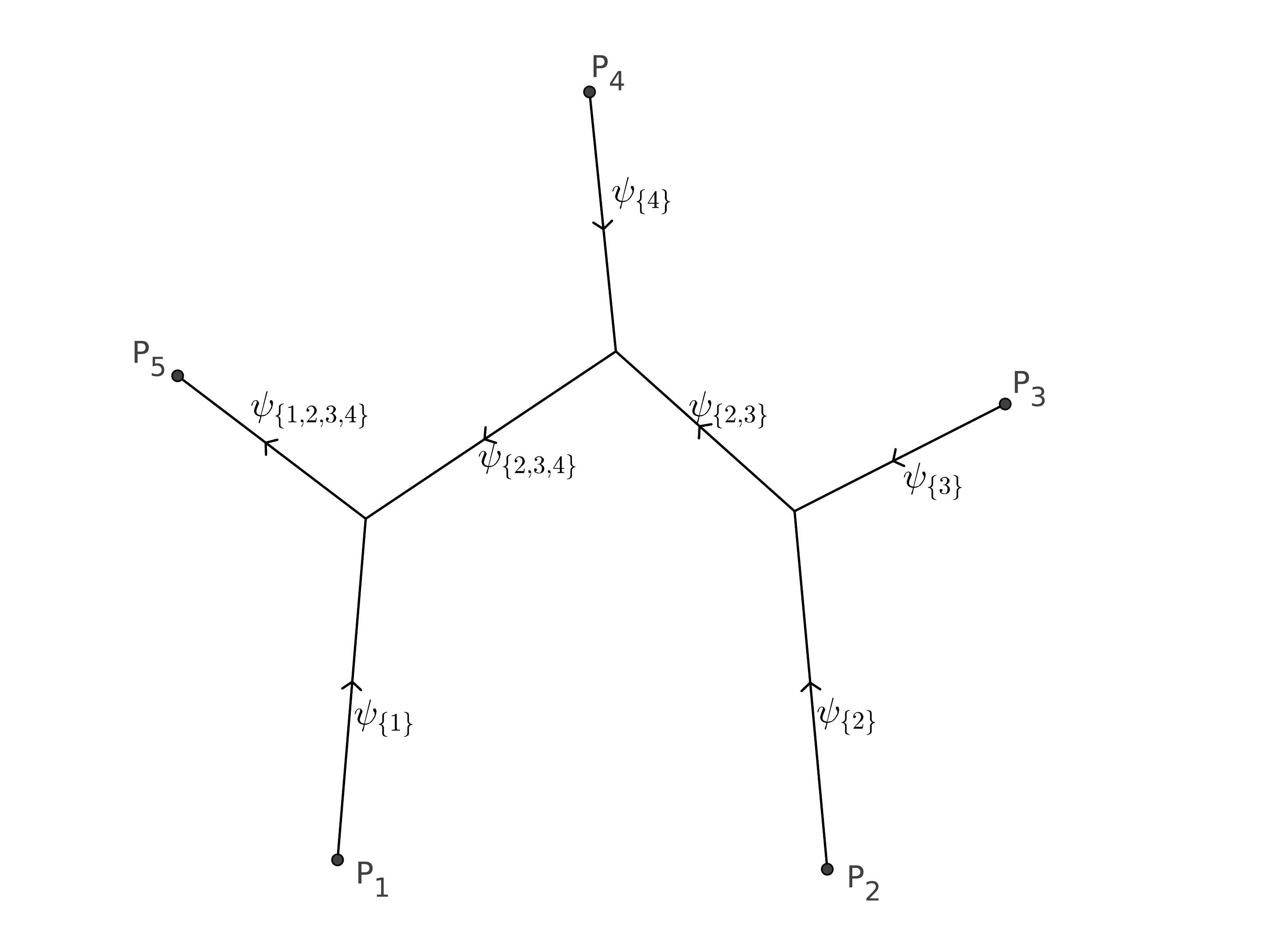}
	\end{tabular}
	\caption{Expected behaviour of the variables $\{\psi_J\}_J$. We can see how each $V_i$ can be reconstructed as the union of the $\psi_J$ such that $i \in J$ and that only a subset of the $\psi_J$ is used.}
	\label{fig:duali}
\end{figure}

Generally speaking, in an optimum for \eqref{eq:psienergy} most of the variables
$\psi_J$ will turn out to be identically $0$ while the ones that are not $0$
everywhere will be concentrated on small regions of the domain. Indeed each
$\psi_J$ can be seen as a possible building block of the final solution because,
due to formula \eqref{eq:Vpsirelation}, the vector field
$\psi_{\{j_1,\dots,j_k\}}$, $\{j_1,\dots,j_k\}\subset \{1,\dots,N-1\}$,
represents the portion of the graph where the fields $V_{j_1}, \dots, V_{j_k}$
coincide (see for example figure \ref{fig:duali} for a visual depiction in two
cases). This means that we expect only a few $\psi_J$ to be non zero on each
element of the grid. With this in mind we can add the following selection
procedure to the previous refinement scheme: given an approximate solution on a
grid $\mathcal{T}$, we identify for each square element $S_m$ which are the non
zero variables $\psi_{J_1^m},\dots,\psi_{J_{k_m}^m}$ on that element and then,
at the next step, we introduce only these variables in that particular region
(in case the element $S_m$ is one of those labelled as ``used'' this means that
in the next optimization we will use only $\psi_{J_1^m},\dots,\psi_{J_{k_m}^m}$
within its $4$ children).

The main advantage of this procedure is clear: once we are able to identify
the regions where each variable $\psi_J$ concentrates (if any) we can
dramatically reduce the number of unknowns we need to introduce, passing from
$2^{N-1}$ vector fields to be defined on each element to only a few of them. Thanks to these two refinement procedures we are now in a position to efficiently tackle the optimization of $\mathcal{R}^\alpha$ using accurate conic solvers.

\section{Results in flat cases}\label{sec:results}

 We present in this section different results obtained using the outlined
scheme integrated with the two refinement procedures described above.

\begin{figure}[tbh]
	\centering
	\begin{tabular}{cccc}
		\includegraphics[width=0.22\textwidth]{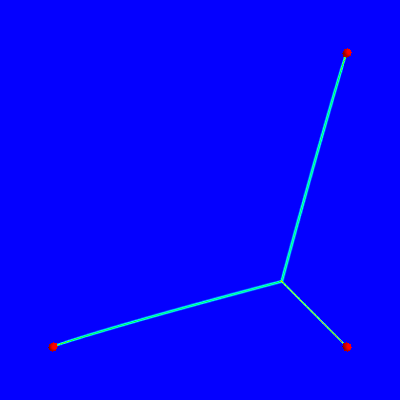} &
		\includegraphics[width=0.22\textwidth]{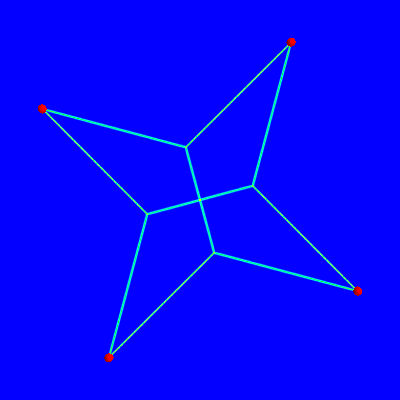} &
		\includegraphics[width=0.22\textwidth]{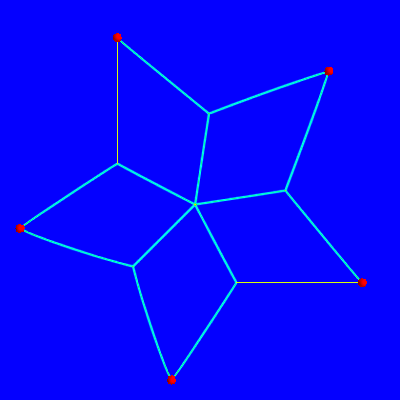} &
		\includegraphics[width=0.22\textwidth]{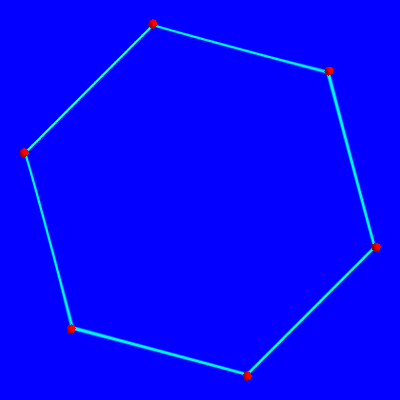}
	\end{tabular}
	\caption{Optima of $\mathcal{R}^0$ in $\R^2$ for $3,4,5,6$ terminal points on the vertices of regular polygons.}
	\label{fig:2d_regular_points}
\end{figure}

In figure \ref{fig:2d_regular_points} we compute minimizers of the relaxed
energy $\mathcal{R}^0$ for regular configurations of terminal points placed on
the vertices of a triangle, a square, a pentagon and an hexagon. In all cases we
start with a regular $32\times 32$ mesh and then apply the previous refinement
procedures $5$ times, ending up with a grid size of $1/1024$ around the optimal
structure. In the first example we are able to retrieve the unique minimizer
while in the second example we obtain a convex combination of the two possible
minimizers for (STP). In the latter case this behaviour is expected because for
this particular configuration of points the relaxation is sharp due to the
calibration argument presented in \cite{MaMa}. In the third experiment we
recover the star-shaped counterexample of figure \ref{fig:steiner_vs_stella}
which seems to be the actual minimizer of the relaxed problem and in the last
picture we get a convex combination of the six possible minimizers. We remark
that the hexagon case is not a calibrated example in the work of
Marchese--Massaccesi but our numerical result suggests the existence of a
calibration because the relaxation seems to be sharp.

\begin{figure}[tbh]
	\centering
	\begin{tabular}{ccc}
		\includegraphics[width=0.22\textwidth]{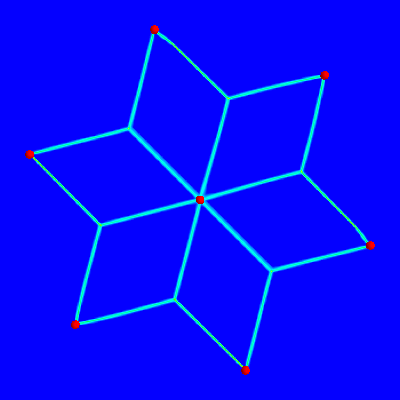} &
		\includegraphics[width=0.22\textwidth]{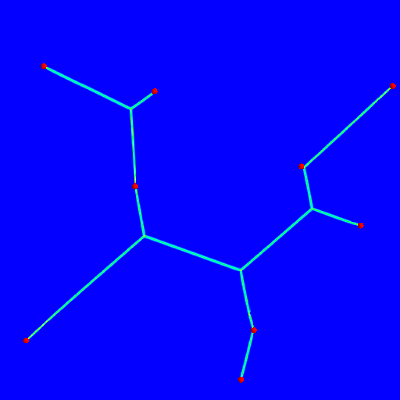} &
		\includegraphics[width=0.22\textwidth]{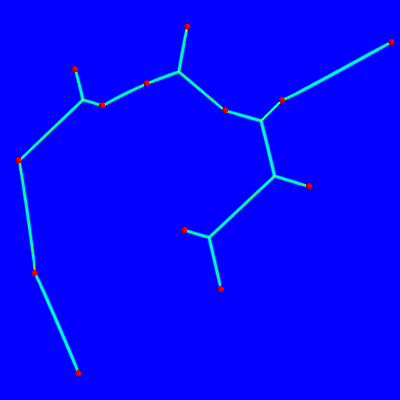}
	\end{tabular}
	\caption{Optima of $\mathcal{R}^0$ in $\R^2$ for $7, 9$ and $13$ terminal points.}
	\label{fig:2d_other_points}
\end{figure}
%
%


 In figure \ref{fig:2d_other_points} we first compute a minimizer for a $7$
points configuration ($6$ vertices of the hexagon plus the center) and observe
how we are able to obtain a convex combination of the two Steiner trees (again
this was expected due to a calibration argument). We observe that in this
example the points do not lie on the boundary of a convex set, meaning that the
problem cannot be simplified into an optimal partition problem as it is done
for example in \cite{ChCrPo}. We then move to some non symmetric distributions
of terminal points: in the second picture we see the result for $9$ randomly
selected points while in the third one we increase the number of terminals up
to $13$. In this last case an ad-hoc approach is necessary. Due to the high
number of variables introduced in \eqref{eq:discrete_energy} a direct
minimization using a conic solver is unfeasible even for very coarse grids
(the amount of memory required to just set up the interior point solver is too
much). To circumvent this problem we first compute a rough solution either
optimizing $\mathcal{R}^0$ on a coarse grid using the primal-dual minimization
scheme or applying the augmented graph idea presented in section \ref{sec:graph}
(see picture \ref{fig:graph_steiner}), and then we use this approximation for deducing
which are the variables $\psi_J$ active at a given point: for every cell
$(k,\ell)$ of a uniform grid we introduce $\psi_{\{j_1,\dots,j_k\}}$ on that
cell only if in the approximate solution every field $V_{j_1},\dots,V_{j_k}$ is
not identically zero in a suitable neighbourhood of the cell. This way we rule
out a huge amount of the $\psi_J$ obtaining a problem which is now tractable
through interior point schemes.

\begin{figure}[tbh]
	\centering
	\begin{tabular}{cccc}
		\includegraphics[width=0.22\textwidth]{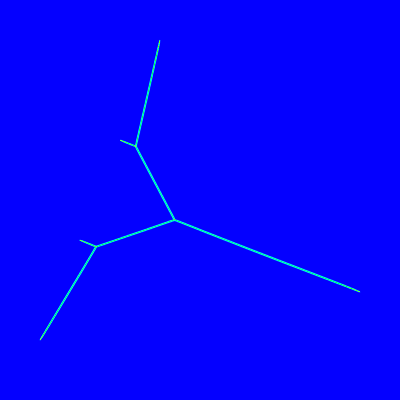} &
		\includegraphics[width=0.22\textwidth]{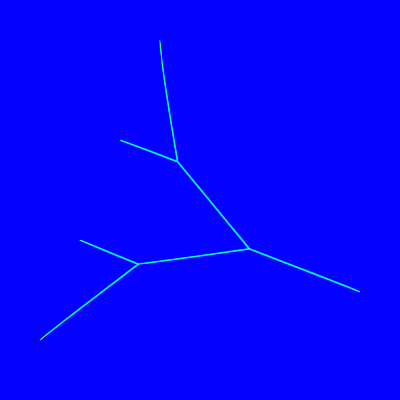} &
		\includegraphics[width=0.22\textwidth]{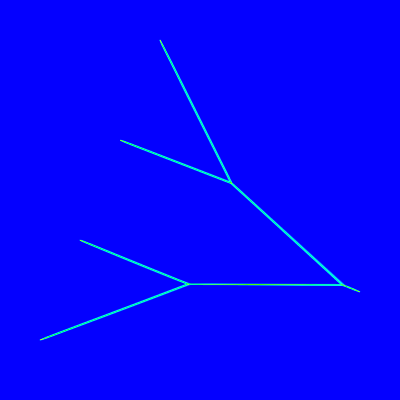} &
		\includegraphics[width=0.22\textwidth]{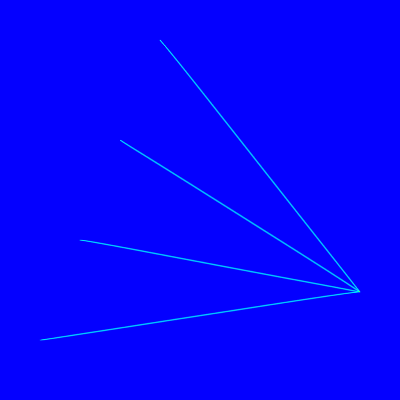}
	\end{tabular}
	\caption{Irrigation networks minimizing $\mathcal{R}^\alpha$ and moving $4$ masses to a unique sink, $\alpha = 0.6, 0.8, 0.95, 1$.}
	\label{fig:2d_irrigation}
\end{figure}

 In figure \ref{fig:2d_irrigation} we test the relaxation $\mathcal{R}^\alpha$
for a simple irrigation problem where we approximate the shape of the optimal
network moving $4$ unit masses located at $S_1=(0.4, 0.9)$, $S_2 = (0.3,0.65)$,
$S_3 = (0.2,0.4)$, $S_4 = (0.1,0.15)$, to the unique sink $T = (0.9,0.27)$. We
can see how for small $\alpha$ the optimal shape is close to the optimal
Steiner tree while for higher values of $\alpha$ the network approaches more
and more the configuration for an optimal Monge--Kantorovitch transport
attaining it for $\alpha = 1$ as expected.

\begin{figure}[tbh]
	\centering
	\begin{tabular}{ccc}
		\includegraphics[width=0.3\textwidth]{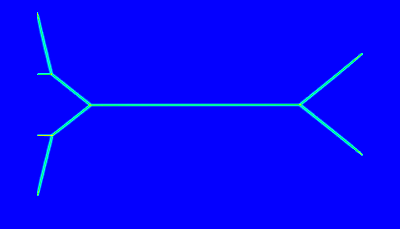} &
		\includegraphics[width=0.3\textwidth]{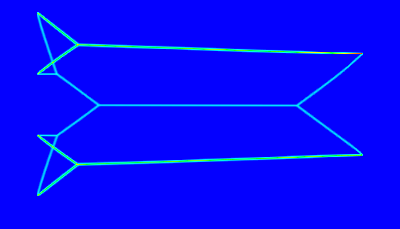} &
		\includegraphics[width=0.3\textwidth]{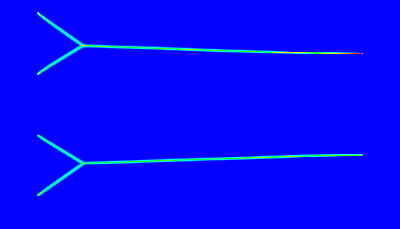}
	\end{tabular}
	\caption{Optima of $\mathcal{R}^\alpha$ for moving $4$ masses from left to right, $\alpha = 0.65, 0.7, 0.75$. The pairings realizing the first infimum and the third one are different.}
	\label{fig:2d_branched}
\end{figure}

 We turn next in figure \ref{fig:2d_branched} to an example where $4$ unit
masses located at $4$ sources on the left ($S_1=(0.1, 0.55)$, $S_2 =
(0.1,0.4)$, $S_3 = (0.1,0.25)$, $S_4 = (0.1,0.1)$) has to be moved to $2$ sinks
of magnitude $2$ on the right ($T_1=(0.9, 0.2)$, $T_2 = (0.9,0.45)$). Since for
each mass we have two possible destinations we need to loop over all feasible
combinations of source/sink couples, solve the corresponding ``who goes where''
problem and then choose the one giving the optimizer with less energy. In the
examples the optimal couplings are
$\{(S_1,T_1),(S_2,T_1),(S_3,T_2),(S_4,T_2)\}$ for $\alpha = 0.65$ and
$\{(S_1,T_2),(S_2,T_2),(S_3,T_1),(S_4,T_1)\}$ for $\alpha = 0.75$. In the case
$\alpha = 0.7$ we are at the switching point between a connected and a
disconnected optimal structure and our relaxed optimum concentrates on both.
\begin{figure}[tbh]
	\centering
	\begin{tabular}{ccc}
		\includegraphics[width=0.23\textwidth]{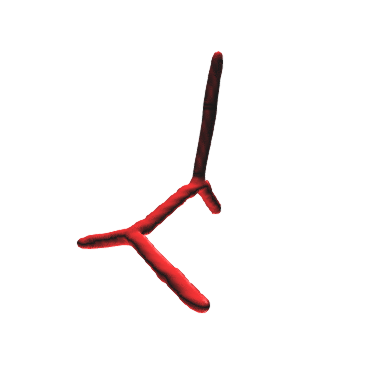} &
		\includegraphics[width=0.23\textwidth]{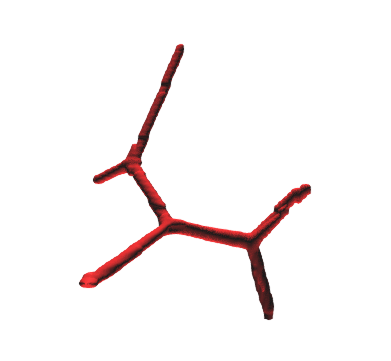} &
		\includegraphics[width=0.23\textwidth]{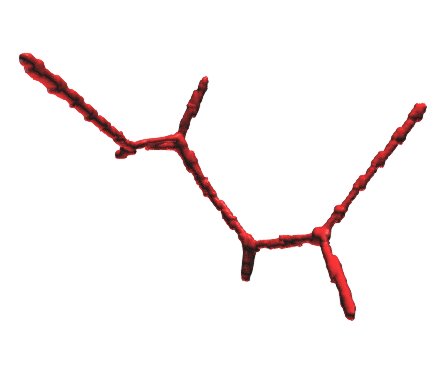}
	\end{tabular}
	\caption{Optima of $\mathcal{R}^0$ for $4,5$ and $7$ points in $\R^3$.}
	\label{fig:3d_steiner}
\end{figure}

 The numerical scheme we have described for the two dimensional case can be
extended directly to the three dimensional context for addressing the
optimization of $\mathcal{R}^\alpha$ in $\R^3$. Non-conformal quadtree type
grids are replaced by non-conformal octree type grids (see \cite{samet88}) and
a staggered approach is employed placing the degrees of freedom on faces of
each cubic element composing the grid. The underlying structure of the discrete
optimization we end up with remains the same and the two refinement procedures
can be extended as they are, without any major change. We see in figure
\ref{fig:3d_steiner} the results for $4,5$ and $7$ points configurations. All
the examples are purely $3$-dimensional and in the first two cases we have a
maximum number of Steiner points (respectively $2$ and $3$), while in the last
case the optimal structure consists of two ``disjoint'' optimal sub-trees
connected through a central terminal point.

\section{Extension to surfaces}\label{sec:surfaces}

As already observed in Section \ref{sec:convexrelaxation} the proposed
relaxation can also be used to address (STP) and $\alpha$-irrigation problems on
surfaces. Up to our knowledge, even in the Steiner tree case, this is the first numerical approximation of these problems covering the manifold framework. Theoretically speaking what we need to do is to solve problem \eqref{eq:relaxed} on a manifold $S$ embedded in $\R^3$, where now a candidate minimizer $\Lambda$ is a matrix valued measure defined on the manifold and divergence constraints translate accordingly. From a numerical point of view our unknowns are again vector fields $(V_1,\dots,V_{N-1})$ living on the surface and the domain will be approximated by means of a triangulated surface $\mathcal{T}_h$. We first discuss the direct extension of the staggered grid idea to $\mathcal{T}_h$ and then present a more accurate discretization, eventually used in our experiments.

\subsection{Raviart--Thomas approach}
The staggered approach presented for quadrilateral grids can be extended to
triangular meshes considering a discretization based on the so-called
Raviart--Thomas basis functions, which are vector valued functions whose degrees
of freedom are related to the flux of the given basis function across edges (see \cite{BrFo}).

 Let $\mathcal{T}_h$ be a regular triangulation of $S$, with $n$ vertices and
$m$ edges, and consider the lowest order Raviart--Thomas basis functions over
$\mathcal{T}_h$: for each edge $e$ in the triangulation we call $K_-$ the
``left'' triangle adjacent to $e$ and $K_+$ the ``right'' triangle adjacent to
$e$ (according to a given fixed orientation) and define the vector function
\[
\Phi_e({x}) =
\left\lbrace
\begin{aligned}
& \frac{\ell_e}{2A_e^+}({x}-{p}_+) & \quad & \text{if } {x} \in K_+ \\
& -\frac{\ell_e}{2A_e^-}({x}-{p}_-) & \quad & \text{if } {x} \in K_- \\
& (0, 0, 0) & \quad & \text{otherwise} \\
\end{aligned}
\right.
\]
where $\ell_e$ is the length of the edge, $A_e^\pm = |K_\pm|$ are the areas of
the triangles and ${p}_+$, ${p}_-$ are the opposite corners (with
the obvious modification for boundary edges). We then approximate each $V_i$, $i =
1,\dots,N-1$, as
\[
V_i({x}) = \sum_{e = 1}^{m} V_i^e \, \Phi_e({x})
\]
and as before matrix valued variables $\phi = (\phi_1,\dots,\phi_{N-1})$ are
considered to be piecewise constant over each element of the triangulation, i.e.
$\phi_i |_K = \phi_i^{K} = (\phi_{i,1}^{K}, \phi_{i,2}^{K}, \phi_{i,3}^{K}) \in
\R^3$ for all $K \in \mathcal{T}_h$, $i = 1,\dots,N-1$. The unknowns are then
the family of parameters $(V_i^e)$ and $(\phi_{i,d}^K)$. Looking at
$\mathcal{R}^\alpha$ the integral we need to compute becomes

\begin{equation}\label{eq:P2discrete}
\sum_{K \in \mathcal{T}_h} \sum_{i=1}^{N-1} \int_{K} \left(\sum_{e = 1}^{m} V_i^e \, \Phi_e({x})\right) \cdot \phi_{i}^K \, d{x},
\end{equation}
and can be made explicit as follows: let $e^K_j$ be the edge of triangle $K$
opposite to point $P_{j}^K$ ($j$-th point of triangle $K$) and $s^{K,e^K_j}=\pm
1$ the position of that triangle with respect to the edge $e_j^K$, then
\eqref{eq:P2discrete} yields
\[
\begin{aligned}
\frac16 \sum_{K \in \mathcal{T}_h} \sum_{i=1}^{N-1} \Big[ & s^{K,e^K_1}\ell_{e^K_1} V_i^{e^K_1}(P_{2}^K+P_{3}^K-2P_{1}^K)\phi^K_{i,1} + s^{K,e^K_2}\ell_{e^K_2} V_i^{e^K_2}(P_{1}^K+P_{3}^K-2P_{2}^K)\phi^K_{i,2} \\ &
+ s^{K,e^K_3}\ell_{e^K_3} V_i^{e^K_3}(P_{1}^K+P_{2}^K-2P_{3}^K) \phi^K_{i,3} \Big].
\end{aligned}
\]
The structure of the discrete energy is the same as the one obtained in the
Euclidean setting (the conditions on $\phi$ translates again in the element-wise
constraint $\phi^K \in K^\alpha$ for all $K\in \mathcal{T}_h$). Furthermore
within this Raviart--Thomas framework we have two advantages: fields $V_i$ are
by construction surface vector fields (i.e. they live in the tangent space to
the surface) and divergence constraints translate into simple flux conditions of
the form
\[
s^{K,e_1^K} \ell_{e_1^K} V_i^{e_1^K} + s^{K,e_2^K} \ell_{e_2^K} V_i^{e_2^K} + s^{K,e_3^K} \ell_{e_3^K} V_i^{e_3^K} = 0 \text{ or } \pm 1
\]
depending on $K$ containing $P_i$, $P_N$ or none of them, and $V_i^e=0$ whenever $e$ is a boundary edge.
The price to pay for such simplicity resides in the fact that this Raviart--Thomas
approximation is a low-order scheme. The objects we would like
to approximate are singular vector fields concentrated on $1$-dimensional
structures but with this approach we generally obtain solutions that are quite
diffused and can only give us an approximate idea of the underlying optimal set.
At the same time this diffusivity prevents a good refinement because the refined
region turns out to be too large. For this reason a better approximation space is needed, even if we will end up with a more complex discrete problem.

\subsection{$\mathbb{P}_2$-based approach}
Let $\mathcal{T}_h$ be a regular triangulation of $S$. We consider the standard discrete space
\[
X_h^2 = \{ v_h \in C^0(\mathcal{T}_h) \,:\, v_h |_K \in \mathbb{P}_2(K),\text{ for all } K \in \mathcal{T}_h \}
\]
and take vector fields $V_i \in (X_h^2)^3$ for all $i = 1,\dots ,N-1$. As in the staggered case matrix valued variables $\phi = (\phi_1,\dots,\phi_{N-1})$ are defined to be piecewise constant over each element of the triangulation, i.e. $\phi_i |_K = \phi_i^{K} =
(\phi_{i,1}^{K}, \phi_{i,2}^{K}, \phi_{i,3}^{K})$ for all $K \in \mathcal{T}_h$, $i = 1,\dots,N-1$. The energy $\mathcal{R}^\alpha$ is then
\[
\sup \left\{ \sum_{K \in \mathcal{T}_h} \sum_{i=1}^{N-1} \int_K V_i \cdot \phi_i^K\,d{x}, \quad \phi^K \in K^\alpha \text{ for all } K \in \mathcal{T}_h\right\}
\]
and the integral over each triangle $K$ can be computed explicitly in terms of
the degrees of freedom associated to $V_i$ and $\phi_i$, $i = 1,\dots,N-1$ (the
integrand reduces to a polynomial of degree $2$). We are left with the
specification of how we impose divergence and tangency constraints on each
$V_i$, $i = 1,\dots,N-1$.

\smallskip

\textbf{Divergence constraints}: for each vector field $V_i$ we have to impose $\div V_i = \delta_{P_i} -
\delta_{P_N}$, where this time the divergence has to be interpreted as the
tangential divergence operator on surfaces (see for instance \cite{Fenics}). We
observe that $\div V_i$ is piecewise linear over each element of the
triangulation and thus, for $K \in \mathcal{T}_h$ not containing $P_i$ or $P_N$
we impose $(\div V_i)|_K = 0$ requiring it to be $0$ at the three vertices of
$K$. On the other hand, if $K\in \mathcal{T}_h$ contains $P_i$ (or $P_N$) we
require the flux of $V_i$ over $\partial K$ to be $+1$ (or $-1$). Eventually,
for each boundary edge $e_b$ of the triangulation we request the flux of $V_i$
through $e_b$ to be $0$.

\smallskip

\textbf{Tangency constraints}: while for the Raviart--Thomas approach the approximate fields are surface
 vector fields by construction, for this $\mathbb{P}_2$ approach we need to
 impose this constraint as an additional condition. For doing so we require
 tangency of $V_i$ at each node of the triangulation and at the mid-point of
 each edge. Normals at these points are approximated as a weighted average of
 normals of surrounding elements.

\medskip
 The above constraints, as it happens in the staggered case, translate into
linear constraints over the degrees of freedom of $V_1,\dots,V_{N-1}$, and the
discrete problem we end up with can be solved using the same strategies
presented in Section \ref{sec:algo}. Eventually we observe that we can extend
the refinement procedures of Section \ref{sec:refinement} also on triangulated
surfaces taking advantage of the re-meshing functionalities of the Mmg Platform
\cite{mmg}: at each step we identify the region where the solution concentrates
the most and then remesh the surface requiring the new mesh to be finer in that
region and coarser elsewhere.


\subsection{Results}

\begin{figure}[tbh]
	\centering
	\begin{tabular}{ccc}
		\includegraphics[width=0.23\textwidth]{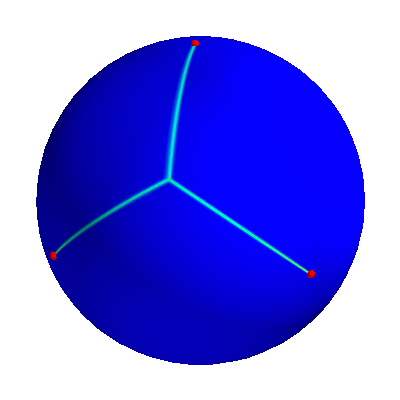} &
		\includegraphics[width=0.23\textwidth]{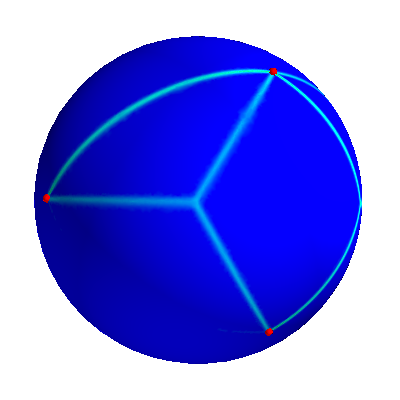} &
		\includegraphics[width=0.23\textwidth]{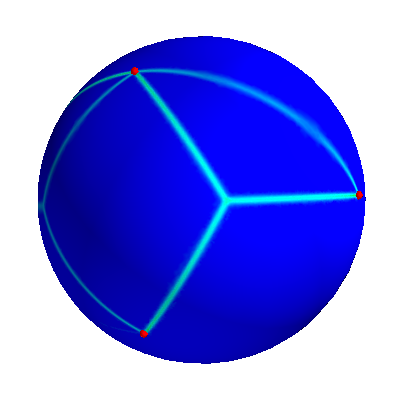}
	\end{tabular}
	\begin{tabular}{cc}
		\includegraphics[width=0.23\textwidth]{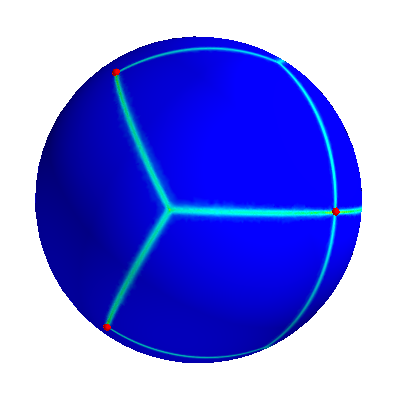} &
		\includegraphics[width=0.23\textwidth]{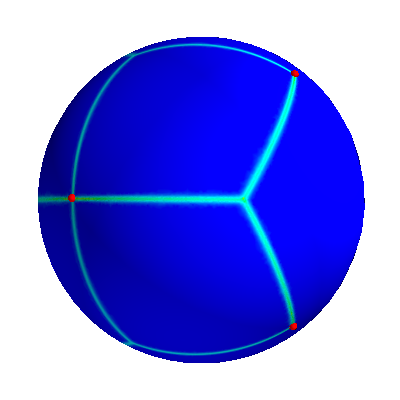}
	\end{tabular}
	\caption{Optima of $\mathcal{R}^0$ for $3,4,5$ points on a sphere, single view for $3$ terminals and different view angles for $4,5$ terminals.}
	\label{fig:3d_sphere}
\end{figure}

In figure \ref{fig:3d_sphere} we see the results obtained through the $\mathbb{P}_2$-based approach for $3$ instances of (STP) on
the sphere. In the first case (upper-left) we approximate the Steiner tree
associated to the terminal points $(1,0,0)$, $(0,1,0)$, $(0,0,1)$, and observe
how we get a classical triple junction. In the second example (upper-middle and
upper-right) we add a fourth point, $(0,-1,0)$, and obtain a convex combination
of minimizers: in this case a possible minimizer can be constructed using the
structure of the first picture completed with an geodesic arc connecting $(0,0,1)$ to
$(0,-1,0)$. We also observe that due to the refinement steps energy concentrates
only on two of the possible four minimizers, the two around which the mesh gets
refined. In the third example (second row) we add a fifth point, $(-1,0,0)$, and
obtain a convex combination of the two minimizers.

\begin{figure}[tbh]
	\centering
	\begin{tabular}{cccc}
		\includegraphics[width=0.21\textwidth]{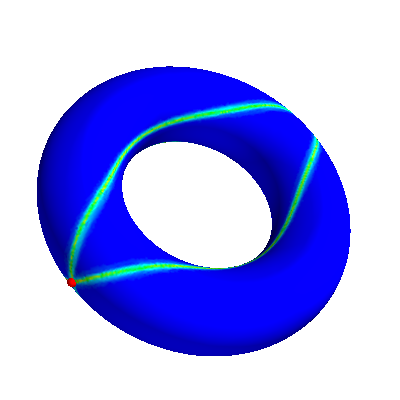} &
		\includegraphics[width=0.21\textwidth]{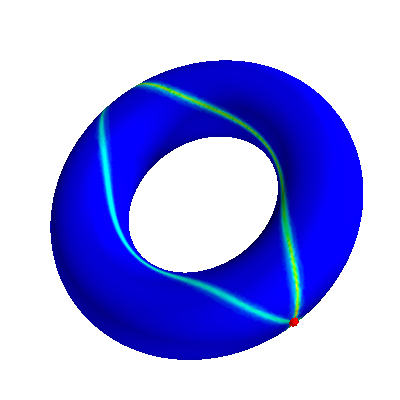} &
		\includegraphics[width=0.21\textwidth]{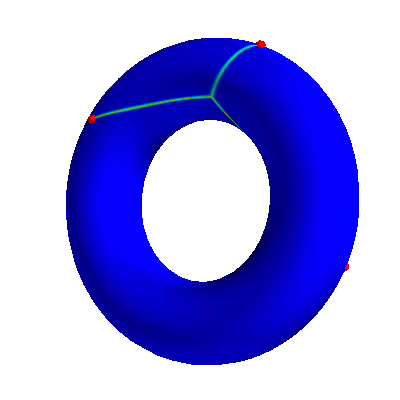} &
		\includegraphics[width=0.21\textwidth]{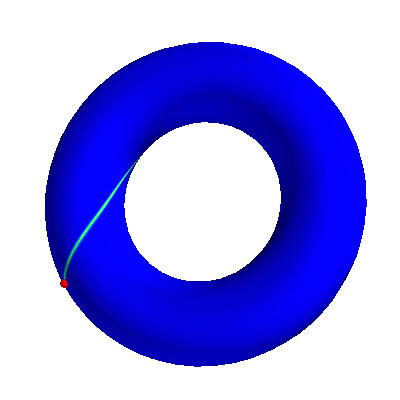}
	\end{tabular}
	\begin{tabular}{ccc}

	\includegraphics[width=0.21\textwidth]{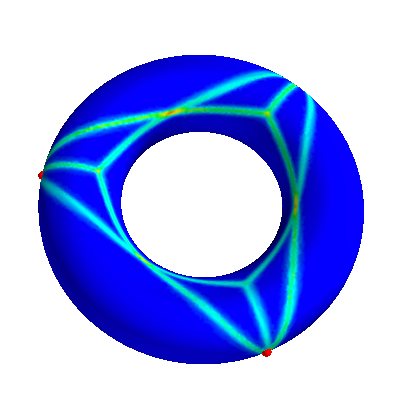} &
	\includegraphics[width=0.21\textwidth]{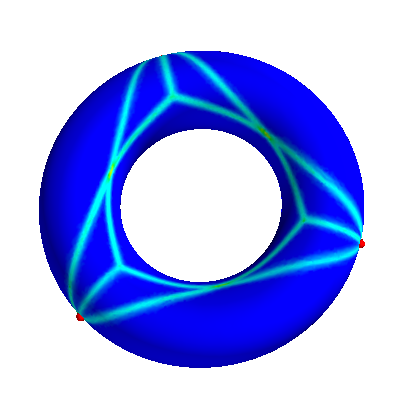} &
	\includegraphics[width=0.21\textwidth]{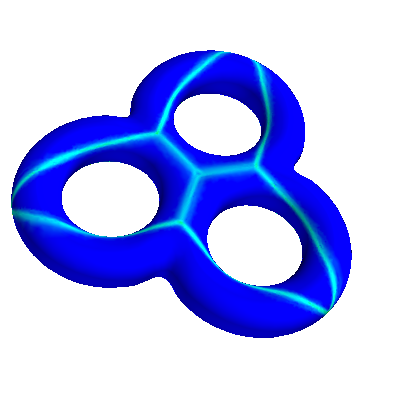}
		\end{tabular}
	\caption{Optima of $\mathcal{R}^0$ for $2,3$ points on different tori (front/back views).}
	\label{fig:3d_tori}
\end{figure}

 As we change the topological nature of the surface results become more
interesting. We approximate in figure \ref{fig:3d_tori} minimizers of
$\mathcal{R}^0$ for some points configurations on the torus. In the first
example (upper-left) we fix two terminal points opposite to each other on the
largest equator and observe an energy concentration on four different paths
(each one a geodesic connecting the two points). For certain $3$ points
configurations we obtain a unique structure with a triple junction
(upper-right), while for $3$ points in a symmetric disposition on the largest
equator we observe as solution a convex combination of the $6$ possible
minimizers (bottom-left). In the last example (bottom-right) we increase the
number of holes of our torus and obtain for a symmetrical $3$ points
configuration a minimizer which cannot be seen as a convex combination of
Steiner trees (i.e. another non sharpness example).

\begin{figure}[tbh]
	\centering
	\begin{tabular}{ccc}
		\includegraphics[width=0.23\textwidth]{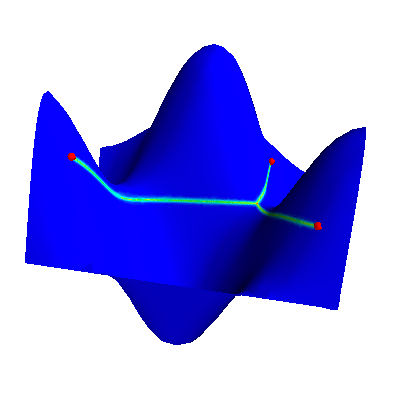} &
		\includegraphics[width=0.23\textwidth]{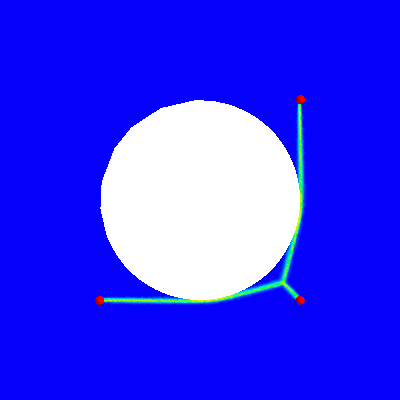} &
		\includegraphics[width=0.23\textwidth]{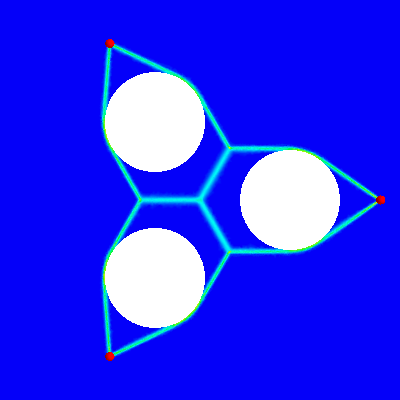}
	\end{tabular}
	\caption{Optima of $\mathcal{R}^0$ for $3$ points on the graph of a function and on some punctured domains in $\R^2$.}
	\label{fig:2d_hloes}
\end{figure}

Finally, in figure \ref{fig:2d_hloes}, we test our relaxation on some surfaces
with boundary. In the first example we connect three given points on the graph
of a function while in the last two we use flat surfaces with holes, which can
be seen as the flat version of the previous tori. In this case solutions can
adhere to the interior boundary of the domain as long as this is energetically
favourable. Observe that, similarly to counter example of figure \ref{fig:steiner_vs_stella}, we obtain a profile which is not a convex combination
 of optimal trees. As in figure \ref{fig:steiner_vs_stella}, we suspect this
  solution to illustrate the fact that our convexification  may be  not sharp in
   specific situations.

\section*{Acknowledgements}

The second author gratefully acknowledges the support of the ANR through the
project GEOMETRYA, the project COMEDIC and  the LabEx PERSYVAL-Lab (ANR-11-LABX-0025-01).



\bibliographystyle{plain}
\bibliography{bibliography}

\begin{thebibliography}{10}

\bibitem{mmg}
{\em Online at https://www.mmgtools.org/}.

\bibitem{AmBr}
Luigi Ambrosio and Andrea Braides.
\newblock Functionals defined on partitions in sets of finite perimeter. {I}.\
  {I}ntegral representation and {$\Gamma$}-convergence.
\newblock {\em J. Math. Pures Appl. (9)}, 69(3):285--305, 1990.

\bibitem{AmBr2}
Luigi Ambrosio and Andrea Braides.
\newblock Functionals defined on partitions in sets of finite perimeter. {II}.\
  {S}emicontinuity, relaxation and homogenization.
\newblock {\em J. Math. Pures Appl. (9)}, 69(3):307--333, 1990.

\bibitem{Ar}
Sanjeev Arora.
\newblock Polynomial time approximation schemes for {E}uclidean traveling
  salesman and other geometric problems.
\newblock {\em J. ACM}, 45(5):753--782, 1998.

\bibitem{Ar1}
Sanjeev Arora.
\newblock Approximation schemes for {NP}-hard geometric optimization problems:
  a survey.
\newblock {\em Math. Program.}, 97(1-2, Ser. B):43--69, 2003.
\newblock ISMP, 2003 (Copenhagen).

\bibitem{bass87}
J.~M. Bass and J.~T. Oden.
\newblock Adaptive finite element methods for a class of evolution problems in
  viscoplasticity.
\newblock {\em Internat. J. Engrg. Sci.}, 25(6):623--653, 1987.

\bibitem{ben2001lectures}
Ahron Ben-Tal and Arkadi Nemirovski.
\newblock {\em Lectures on modern convex optimization: analysis, algorithms,
  and engineering applications}, volume~2.
\newblock Siam, 2001.

\bibitem{BeCaMo}
Marc Bernot, Vicent Caselles, and Jean-Michel Morel.
\newblock {\em Optimal transportation networks: models and theory}, volume
  1955.
\newblock Springer Science \& Business Media, 2009.

\bibitem{BoOrOu}
Mauro Bonafini, Giandomenico Orlandi, and \'Edouard Oudet.
\newblock Variational approximation of functionals defined on $1$-dimensional
  connected sets: the planar case.
\newblock {\em SIAM J. Math. Anal.}, accepted.

\bibitem{BoBrLe}
Matthieu Bonnivard, Elie Bretin, and Antoine Lemenant.
\newblock Numerical approximation of the steiner problem in dimension 2 and 3.
\newblock 2018.

\bibitem{BoLeSa}
Matthieu Bonnivard, Antoine Lemenant, and Filippo Santambrogio.
\newblock Approximation of length minimization problems among compact connected
  sets.
\newblock {\em SIAM J. Math. Anal.}, 47(2):1489--1529, 2015.

\bibitem{BoVa}
Guy Bouchitt{\'e} and Michel Valadier.
\newblock Integral representation of convex functionals on a space of measures.
\newblock {\em Journal of functional analysis}, 80(2):398--420, 1988.

\bibitem{BrFo}
Franco Brezzi and Michel Fortin.
\newblock {\em Mixed and hybrid finite element methods}, volume~15.
\newblock Springer Science \& Business Media, 2012.

\bibitem{ChCrPo}
Antonin Chambolle, Daniel Cremers, and Thomas Pock.
\newblock A convex approach to minimal partitions.
\newblock {\em SIAM J. Imaging Sci.}, 5(4):1113--1158, 2012.

\bibitem{ChMeFe}
Antonin Chambolle, Luca Alberto~Davide Ferrari, and Benoit Merlet.
\newblock A phase-field approximation of the steiner problem in dimension two.
\newblock {\em Advances in Calculus of Variations}, 2017.

\bibitem{ChMeFe2}
Antonin Chambolle, Luca Alberto~Davide Ferrari, and Benoit Merlet.
\newblock Variational approximation of size-mass energies for k-dimensional
  currents.
\newblock {\em arXiv preprint arXiv:1710.08808}, 2017.

\bibitem{JuMP}
Iain Dunning, Joey Huchette, and Miles Lubin.
\newblock Ju{M}{P}: A {M}odeling {L}anguage for {M}athematical {O}ptimization.
\newblock {\em SIAM Review}, 59(2):295--320, 2017.

\bibitem{dykstra83}
Richard~L Dykstra.
\newblock An algorithm for restricted least squares regression.
\newblock {\em Journal of the American Statistical Association},
  78(384):837--842, 1983.

\bibitem{MINLP}
Claudia D’Ambrosio, Marcia Fampa, Jon Lee, and Stefan Vigerske.
\newblock On a nonconvex minlp formulation of the euclidean steiner tree
  problem in n-space.
\newblock In {\em International Symposium on Experimental Algorithms}, pages
  122--133. Springer, 2015.

\bibitem{Gi}
Edgar~N Gilbert.
\newblock Minimum cost communication networks.
\newblock {\em Bell Labs Technical Journal}, 46(9):2209--2227, 1967.

\bibitem{Karp}
Richard~M Karp.
\newblock Reducibility among combinatorial problems.
\newblock In {\em Complexity of computer computations}, pages 85--103.
  Springer, 1972.

\bibitem{MaMa2}
Andrea Marchese and Annalisa Massaccesi.
\newblock An optimal irrigation network with infinitely many branching points.
\newblock {\em ESAIM Control Optim. Calc. Var.}, 22(2):543--561, 2016.

\bibitem{MaMa}
Andrea Marchese and Annalisa Massaccesi.
\newblock The {S}teiner tree problem revisited through rectifiable
  {$G$}-currents.
\newblock {\em Adv. Calc. Var.}, 9(1):19--39, 2016.

\bibitem{MaOuVe}
Annalisa Massaccesi, \'Edouard Oudet, and Bozhidar Velichkov.
\newblock Numerical calibration of {S}teiner trees.
\newblock {\em Applied Mathematics \& Optimization}, pages 1--18, 2017.

\bibitem{mosek2010mosek}
APS Mosek.
\newblock The {MOSEK} optimization software.
\newblock {\em Online at http://www. mosek. com}, 54, 2010.

\bibitem{OuSa}
Edouard Oudet and Filippo Santambrogio.
\newblock A {M}odica-{M}ortola approximation for branched transport and
  applications.
\newblock {\em Arch. Ration. Mech. Anal.}, 201(1):115--142, 2011.

\bibitem{ChPo}
Thomas Pock and Antonin Chambolle.
\newblock Diagonal preconditioning for first order primal-dual algorithms in
  convex optimization.
\newblock In {\em Computer Vision (ICCV), 2011 IEEE International Conference
  on}, pages 1762--1769. IEEE, 2011.

\bibitem{Fenics}
Marie~E. Rognes, David~A. Ham, Colin~J. Cotter, and Andrew~T.T. McRae.
\newblock Automating the solution of {P}{D}{E}s on the sphere and other
  manifolds in {F}{E}ni{C}{S} 1.2.
\newblock {\em Geoscientific Model Development}, 6(6):2099--2119, 2013.

\bibitem{samet88}
Hanan Samet.
\newblock An overview of quadtrees, octrees, and related hierarchical data
  structures.
\newblock In {\em Theoretical Foundations of Computer Graphics and CAD}, pages
  51--68. Springer, 1988.

\bibitem{ScLeBa}
Mark Schmidt, Nicolas~L Roux, and Francis~R Bach.
\newblock Convergence rates of inexact proximal-gradient methods for convex
  optimization.
\newblock In {\em Advances in neural information processing systems}, pages
  1458--1466, 2011.

\bibitem{warme2001geosteiner}
DM~Warme, Pawel Winter, and Martin Zachariasen.
\newblock Geo{S}teiner 3.1.
\newblock {\em Department of Computer Science, University of Copenhagen
  (DIKU)}, 2001.

\bibitem{Xia}
Qinglan Xia.
\newblock Optimal paths related to transport problems.
\newblock {\em Commun. Contemp. Math.}, 5(2):251--279, 2003.

\end{thebibliography}

\end{document}